\documentclass[preprint]{elsarticle}
\journal{Applied Numerical Mathematics}

\usepackage[latin1]{inputenc}
\usepackage{amsmath}
\usepackage{amsfonts}
\usepackage{amssymb}
\usepackage{amsthm}
\usepackage{graphicx}
\usepackage{color}
\usepackage[english]{babel}
\usepackage{hyperref}

\newtheorem{corollary}{Corollary}[section]
\newtheorem{theorem}{Theorem}[section]
\newtheorem{proposition}{Proposition}[section]

\newcommand{\R}{{\ensuremath{\mathbb{R}}}}

\newcommand{\bw}{\boldsymbol w}
\newcommand{\bW}{\boldsymbol W}
\newcommand{\bc}{\boldsymbol c}
\newcommand{\bC}{\boldsymbol C}

\newcommand{\bq}{\boldsymbol q}
\newcommand{\bigO}{\mathcal{O}}

\def\mig{\frac12}
\def\boc#1{\boldsymbol{{\cal #1}}}
\def\mD{\boc{D}}
\def\XXX#1#2{\chi(\cdot, #2)^{-1}(#1)}
\def\ffc#1#2{\frac{\partial\chi_{#1}}{\partial \xi}(c_{#1}, #2)}
\def\ffp#1#2{\frac{\partial\chi_{#1}}{\partial p}(c_{#1}, #2)}
\def\sii{\Leftrightarrow}

\begin{document}

The version of record of this article, first published in Applied Numerical Mathematics, is available online at Publisher's website:
\medskip

\url{https://doi.org/10.1016/j.apnum.2024.03.006}

\newpage
\begin{frontmatter}

\title{WENO scheme on characteristics for the equilibrium dispersive
  model of chromatography with generalized Langmuir isotherms.}
\author{R. Donat \footnote{Email: donat@uv.es}, M.C. Mart\'i \footnote{Email: Maria.C.Marti@uv.es} and P. Mulet \footnote{Email: mulet@uv.es}\\  {\small Department of Mathematics.
Universitat de Val\`encia} \\  {\small Av. Vicent Andr\'es Estell\'es
s/n; 46100 Burjassot (Val\`encia), Spain}}

\begin{abstract}

Column chromatography is a laboratory and industrial technique used to separate different substances mixed in a solution. Mathematically, it can be modelled using non-linear partial
differential equations whose main ingredients are the {\it adsorption
  isotherms}, which are non-linear functions modelling the affinity between the different substances in the solution and the solid stationary phase filling the column.
The goal of this work is twofold. Firstly, we aim to extend the
techniques of Donat, Guerrero and Mulet (Appl. Numer. Math. 123 (2018)
22-42) to other adsorption isotherms. In particular, we propose a
family of generalized {\it Langmuir}-type isotherms and prove that the
correspondence between the concentrations of solutes in the liquid phase (the primitive
variables) and the conserved variables is well defined and admits a
global smooth inverse that can be computed numerically.  Secondly, to
establish the well-posedness of the mathematical model, we study the
eigenstructure of the Jacobian of the mentioned correspondence and
use this characteristic information to get oscillation-free
sharp interfaces on the numerical approximate solutions. To do so, we
determine the structure of the Jacobian matrix of the system and use
it to deduce its eigenstructure. We combine the use of
characteristic-based numerical fluxes with a second-order implicit-explicit
scheme proposed in the cited reference and perform some numerical
experiments with T\'oth's adsorption isotherms to demonstrate that the characteristic-based schemes produce accurate numerical solutions with no oscillations, even when steep gradients appear in the solutions.

\end{abstract}

\end{frontmatter}

\section{Introduction}
Separation processes of complex mixtures are often modelled by non-linear partial differential equations (PDEs). In liquid batch chromatography, the interaction of the {\em solute} (a fluid mixture) with a porous medium that fills a long column results in a separation process that leads, for sufficiently long columns, to band profiles of single-component fluid. The separation process in the mixture is governed by the strength of the interaction of the fluid components with the solid (stationary) phase. Liquid batch chromatography is required when highly pure components need to be extracted from a fluid mixture, for example, in the pharmaceutical industry.

Several mathematical models for the analysis and numerical simulation of liquid batch chromatography can be found in the literature, see e.g. \cite{Mazzotti2013, Guiochon}. Under reasonable assumptions, these models involve a system of convection-dominated PDEs, coupled by the so-called {\em adsorption isotherms}, non-linear functions that establish equilibrium relationships between the solid and liquid phase concentrations. The non-linear character of these functions hinders finding analytical solutions to chromatographic models, even in the simplest situations, and toughens the analysis of their properties. This makes the development of robust numerical tools necessary to obtain accurate numerical solutions of the chromatographic models. The development of such reliable numerical techniques can be undeniably useful in real-life chemistry laboratories, for instance, to reproduce dynamical scenarios including the formation of band profiles of pure components during the separation process without performing trial-and-error empirical experimentation, which is more costly.

In \cite{DGM18}, the authors considered the Equilibrium Dispersive
(ED) model with multicomponent{\em Langmuir} adsorption isotherms. The ED model, which is described in detail in \cite{Guiochon}, is based on the assumption that the mobile and solid phases are in permanent equilibrium at all positions in the column. This equilibrium is described by  {\em adsorption isotherms} $\bq=\bq(\bc)$, where ${\boldsymbol c}=(c_1,\ldots,c_N)^T$ and ${\boldsymbol q}=(q_1,\ldots,q_N)^T$ are the vectors with components $c_i$ and $q_i$ defining the concentration of the $i$-th component of the solution in the liquid and solid phase, respectively.

As observed in \cite{DGM18}, the mass balance equations of the ED model are given by
\begin{align}\label{eq:modelED}
  &\frac{\partial {\boldsymbol c}}{\partial
  t}+\frac{1-\epsilon}{\epsilon}\frac{\partial {\boldsymbol
    q}}{\partial t}+u\frac{\partial {\boldsymbol c}}{\partial z}=D_a
    \frac{\partial^2 {\boldsymbol c}}{\partial^2 z},
\end{align}
which  can be rewritten as
\begin{align}
\label{eq:modelwc}
 &\frac{\partial \bw}{\partial t} + \frac{\partial (u\bc)}{\partial z}=
D_a \frac{\partial^2 \bc}{\partial^2 z}, \quad \bw=\bW(\bc)=\bc+\frac{1-\epsilon}{\epsilon} \bq(\bc),
\end{align}
where $u$ is the velocity of the mobile phase, $\epsilon$ is the total porosity of the solid phase, $0<\epsilon \le 1$, $D_a$ is the {\em apparent axial dispersion coefficient} (under the assumption that band broadening is only caused by axial dispersion), $t$ is time and $z$ denotes the normalised position in the column, where the top corresponds to $z=0$ and the bottom to $z=1$.

The parameter $D_a$ is often very small or even zero, hence sharp discontinuous fronts may appear in the numerical solutions of the non-linear system (\ref{eq:modelED}) that will travel through the column. This makes necessary the use of reliable and efficient numerical techniques that can accurately recover these fronts and also capture their speed of propagation correctly. The existence of a smooth, globally well-defined, one-to-one correspondence between $\bw$ and $\bc$ for {\em Langmuir}  adsorption isotherms was shown in \cite{DGM18}. This allowed the authors to propose and implement conservative numerical schemes (see \cite{Levequeb}) for \eqref{eq:modelwc}, for which the mass conservation of the {\em conserved variables} $\bw$ is guaranteed and the correct speed of propagation of steep fronts is obtained as well.

In this paper,  we propose multicomponent adsorption isotherms which are generalizations of the multicomponent {\em Langmuir} isotherms and prove that they give well-posed initial boundary-value problems.  Specifically, the adsorption isotherms that we propose read as:
\begin{equation} \label{eq:q-def}
q_i(\bc)= \frac{a_i c_i}{\varphi(\sum_{i=1}^N b_i c_i)}, \qquad i=1,2,\dots,N,
\end{equation}
where $\varphi$ is a suitable real function, whose mathematical properties will be stated below, $a_i=\alpha_i b_i$, $\alpha_i>0$ is the column saturation capacity of component $i$ and $b_i>0$ is the ratio of the rate constants of adsorption and desorption for component $i$.
As an instance, T\'oth's isotherms (see \cite{Guiochon}) is given by \eqref{eq:q-def} for $\varphi(c)=\big(1+c^{\nu}\big)^{\frac{1}{\nu}}$, where $0<\nu\leq 1$ is the heterogeneity parameter. Note that {\it Langmuir} isotherms correspond to $\nu=1$.

The main goal of this paper is to generalize the results in \cite{DGM18} for these adsorption isotherms. In this direction, we will begin by showing that, also in this case, $\bW$ is a bijection between the vector of concentrations $\bc$ and the vector of  {\em conserved variables} $\bw$, so that we can rewrite \eqref{eq:modelwc} as:
\begin{align} \label{eq:modelwc1}
&\frac{\partial \bw}{\partial t} + \frac{\partial \boldsymbol{f}(\bw)}{\partial z}=D_a \frac{\partial^2 \bC(\bw)}{\partial^2 z}, \qquad \boldsymbol{f}(\bw)=u \bC(\bw),\end{align}
where $\bC=\bW^{-1}$ is a continuously differentiable function whose explicit expression cannot be determined for $N>1$, but can be approximated efficiently from the only positive root of a certain rational function.

To determine the well-posedness of the system it is necessary to study the eigenvalues of   $\bW'(\bc)$, which are proportional to the inverses of the eigenvalues of the Jacobian matrix of the convective fluxes $u\bC'(\bw)$ and of the diffusion matrix $D_a\bC'(\bw)$ of system \eqref{eq:modelwc1}.  Following \cite{Donat3}, we will prove that the eigenvalues of $\bW'(\bc)$  are strictly positive and pairwise different, and hence the model is strictly hyperbolic when $D_a=0$ and parabolic otherwise.

These theoretical results constitute the key to achieving the second goal of this paper, which is to incorporate characteristic information in the convective numerical fluxes of the fully conservative numerical schemes that we propose to simulate the ED model with generalized Langmuir adsorption isotherms. The characteristic-based discretization of the convective flux uses the spectral decomposition of the Jacobian matrix of the flux to compute the numerical approximations by local projections onto characteristic fields. The numerical solutions obtained with characteristic-based Weighted Essentially Non-Oscillatory (WENO) reconstructions, combined with a correct upwinding strategy, are known to be more precise, in terms of resolution and oscillatory behavior, than the ones obtained using schemes that do not use characteristic projections. To illustrate this fact, in this work we compare the results obtained with the characteristic-based WENO scheme with the component-wise one proposed in \cite{DGM18}, along with other options.

The paper is organized as follows: Section \ref{sec:Modeleq} is devoted to proposing new adsorption isotherms that generalize multicomponent Langmuir isotherms and to the mathematical analysis of the ED  model with these isotherms.  The construction of the characteristic-based WENO schemes used in this work is detailed in Section \ref{sec:charweno}.
In Section \ref{sec:imex-weno}, we describe the second-order implicit-explicit (IMEX) scheme and the different discretizations of the convective terms that will be used in this work, discussing several details for their implementation in numerical simulations of chromatographic processes that fit into the ED model \eqref{eq:q-def}-\eqref{eq:modelwc1}. Some numerical experiments to test the performance of the characteristic-based scheme are shown in Section \ref{sec:numex}. Finally, some conclusions and perspectives for future work are drawn in Section \ref{sec:conc}.

\section{The mathematical structure of the ED model.}\label{sec:Modeleq}

We will start this work analyzing the {\em conservative formulation} of the  ED model \eqref{eq:q-def}-\eqref{eq:modelwc1}, with $\bw= \bW(\bc)= \bc + \frac{1-\varepsilon}{\varepsilon} \bq(\bc)$. The components of $\bW$ are \begin{equation} \label{eq:Wi-defToth}
W_i (\bc)= c_i \left( 1+\frac{\eta_i}{\varphi(\sum_{j}b_j c_j)}
\right),\quad
 \eta_i=\frac{1-\epsilon}{\epsilon}a_i, \quad 1 \leq i\leq N,
\end{equation}
which have the following algebraic structure
\begin{equation}\label{eq:report2-33}
\left\{\begin{aligned}
    &W_i(\bc)=\chi_i(c_i, \varphi(b^T\bc)), \\
    &\chi_i(\xi, p)=\xi \left( 1+\frac{\eta_i} {p} \right),  \quad \chi_i:
  [0, \infty) \times (0,\infty)  \rightarrow [0,\infty),\\
  &\varphi \colon [0,\infty)\to [1,\infty).
 \end{aligned}\right.
\end{equation}

In the rest of the paper, we assume the following hypothesis on the function $\varphi$:
\begin{enumerate}
\item $\varphi$  is a continuous bijection,
\item $\varphi'$ exists and is continuous in $(0,\infty)$ and
  \begin{equation}\label{eq:cphi}
    \varphi'(c) >0,\quad
    \left(\frac{c}{\varphi(c)}\right)'=\frac{\varphi(c)-c\varphi'(c)}{\varphi(c)^2} > 0,\,\forall c\in(0,\infty).
  \end{equation}
\end{enumerate}
In particular, $\varphi(0)=1$ and $\varphi(c)>1,\,\forall c>0$.

T\'oth's isotherms, which fit into this framework with $\varphi(c)=\big(1+c^{\nu}\big)^{\frac{1}{\nu}}$, satisfy these requirements for $\nu>0$, since
\begin{align*}
\varphi'(c)&=\frac{\big(1+c^{\nu}\big)^{\frac{1}{\nu}-1}}{c^{1-\nu}}
  > 0,\,\,\forall c>0,\\
\left(\frac{c}{\varphi(c)}\right)'&=
  \frac{1}{(1+c^{\nu})^{\frac{1}{\nu}+1}} > 0,\,\,\forall c>0.
  \end{align*}

The following general analysis applies to the equations of the ED
model \eqref{eq:q-def}-\eqref{eq:modelwc1} with the notation in \eqref{eq:report2-33}.

\begin{theorem}\label{th:winvertible}
Let $\mathcal{I}, \mathcal{J}, \mathcal{P}$ be real intervals and
\begin{align*}
\chi&\colon \mathcal{I}^N\times \mathcal{P} \to   \mathcal{J}^N,\\
\psi&\colon \mathcal{I}^N\to \mathcal{P},  \\
\bW&\colon \mathcal{I}^N \to \mathcal{J}^N,\quad
\bW(\bc)=\chi(\bc, \psi(\bc)).
\end{align*}
Assume that for any $p\in\mathcal{P}$ the function $\chi(\cdot, p)\colon \mathcal{I}^N\to\mathcal{J}^N$ is bijective and that for any $\bw\in\mathcal{J}^N$ there exists a unique $p=p(\bw)\in\mathcal{P}$ such that
\begin{align}\label{eq:report2-4}
 p&=\psi\Big(\chi\big(\cdot, p\big)^{-1}(\bw)\Big).
\end{align}
Then $\bW$ is bijective and  the function $ \bC\colon\mathcal{J}^N\to\mathcal{I}^N, $ defined as
\begin{equation}\label{eq:cw-def}
\bC(\bw)=\chi(\cdot, p(\bw))^{-1}(\bw),
\end{equation}
is the inverse of $\bW$.
\end{theorem}

\begin{proof}
The function $\bC: \mathcal{J}^N \rightarrow \mathcal{I}^N$ in
\eqref{eq:cw-def} is well defined since the function $\chi(\cdot,
p(\bw))\colon \mathcal{I}^N\to\mathcal{J}^N$ is a bijection by hypothesis.

From \eqref{eq:report2-4} and \eqref{eq:cw-def} we have
\begin{align*}
\psi(\bC(\bw))=\psi(\XXX{\bw}{p(\bw)})=p(\bw), \quad  \chi(\bC(\bw), p(\bw))=\bw,
\end{align*}
therefore
\begin{align*}
  \bW(\bC(\bw))&=\chi(\bC(\bw), \psi(\bC(\bw)))=\chi(\bC(\bw), p(\bw))=\bw.
\end{align*}

Now, given $\bc \in \mathcal{I}^N$, since $\bW(\bc)=\chi(\bc, \psi(\bc))$ and  $\psi(\bc) \in \mathcal{P}$, hence $\chi(\cdot,\psi(\bc))$ is invertible, we have that
\begin{align*}
  \XXX{\bW(\bc)}{\psi(\bc)}=\bc,
\end{align*}
that is $\psi(\bc) =  \psi(\XXX{\bW(\bc)}{\psi(\bc)}), \,\, \forall \bc \in \mathcal{I}^N$, therefore
\begin{align*}
p(\bW(\bc))=\psi(\bc),
\end{align*}
which leads to
\begin{align*}
\bC(\bW(\bc))&=\XXX{\bW(\bc)}{p(\bW(\bc))}=   \XXX{\chi(\cdot, \psi(\bc))(\bc)}{\psi(\bc)}=\bc.
\end{align*}
\end{proof}

As observed in \cite{DGM18} for the multicomponent Langmuir isotherms, the  function
${\boldsymbol C}({\boldsymbol  w})$ in \eqref{eq:cw-def} cannot be
explicitly determined when  $N>1$.  However, the previous result
states that  the value of $\bC(\bw)$ can be efficiently
approximated in practice, for each $\bw \in \mathcal{J}^N$, from the
determination of the unique $p(\bw) \in {\cal P}$ satisfying
\eqref{eq:report2-4}.  For {\it Langmuir} adsorption isotherms, the
computation of this value can be carried out (numerically)  by using
a root finder on a rational function \cite{DGM18}. We shall see that
the determination of $\bC(\bw)$ can also be carried out efficiently in
the general case.

\begin{corollary}
The function $\bW\colon \mathcal{I}^N\to\mathcal{I}^N$, whose components are given in \eqref{eq:Wi-defToth}, is bijective for $\mathcal{I}$ either $(0,\infty) $ or $[0,\infty)$, with inverse given by  \eqref{eq:cw-def}.
\end{corollary}

\begin{proof}
We aim to apply  Theorem \ref{th:winvertible}
with $\mathcal{J}=\mathcal{I},\,\, \mathcal{P}=[1,\infty)$, $\psi(\bc)=\varphi(b^T \bc)$. The function $\chi(\cdot, p)$ from \eqref{eq:report2-33}  is  invertible for any $p>0$, with
\[ \big((\chi(\cdot,p))^{-1}(\bw)\big)_i = \frac{p}{p+\eta_i} w_i, \quad \forall \bw \in \mathcal{I}^N. \]

We check now that $\forall \bw \in \mathcal{I}^N$,  there exists a unique $p \in [1,\infty)$ such that $p= \psi \big(\chi(\cdot, p)^{-1}(\bw)\big)$, which in this context reads:
\begin{align}\label{eq:report2-34}
  p=\varphi\big(\sum_{i=1}^{N}b_iw_i\frac{p}{p+\eta_i}\big)\,\,\sii\,\,
  0=\sum_{i=1}^{N}\frac{b_iw_i}{p+\eta_i}-\frac{\varphi^{-1}(p)}{p}=S_{\bw}(p),
\end{align}
where $S_{\bw}(p)$ is a continuous and strictly decreasing function, since, for $p\in(1,\infty),  c=\varphi^{-1}(p)\in (0,\infty)$,  \eqref{eq:cphi} and $b_iw_i\geq 0\,\forall i,$ imply:
\begin{align*}
S_{\bw}'(p)  =-\sum_{i=1}^{N}\frac{b_iw_i}{(p+\eta_i)^2}-\frac{\varphi(c)-c\varphi'(c)}{\varphi(c)^2\varphi'(c)}<0.
\end{align*}

With the notation
\begin{equation}\label{eq:pb}
  \overline{p}(\bw)=\varphi(b^T\bw),
\end{equation}
we  have
\begin{align*}
&S_{\bw}(\overline{p}(\bw)) \leq \sum_{i=1}^{N}\frac{b_iw_i}{\overline{p}(\bw)}-
   \frac{\varphi^{-1}(\varphi(b^T\bw))}{\overline{p}(\bw)}=0.
\end{align*}
Since  $\varphi(0)=1$ we deduce:
\begin{align*}
&S_{\bw}(1)=\sum_{i=1}^{N}\frac{b_iw_i}{1+\eta_i}\geq 0,
\end{align*}
hence, there exists a unique $p\in[1,\overline{p}(\bw)]$ solving
\eqref{eq:report2-34} by the continuity and strict monotonicity of $S_{\bw}(p)$.

Theorem \ref{th:winvertible} thus applies, and the function
\[ \bC(\bw)= (\chi(\cdot, p(\bw)))^{-1} (\bw), \,\,  \bw \in \mathcal{I}^N, \]
satisfies $\bC=\bW^{-1}$.
\end{proof}

The practical computation of the inverse $\bC(\bw)$ only requires finding $p(\bw)$, the only root of \eqref{eq:report2-34} in $[1, \overline{p}(\bw)]$, with $\overline{p}(\bw)$ given by \eqref{eq:pb}, which can be efficiently found by using a combination of bisection and Newton's method.

The  {\em structure} of the Jacobian matrix of the function $\bW(\bc)= \chi(\bc, \varphi(b^T\bc))$ can also be determined for the general framework stated in Theorem \ref{th:winvertible}, if we further assume, as is the case for \eqref{eq:report2-33}, that $\big(\chi(\bc, p)\big)_{i}=\chi_i(c_i, p)$, for continuously differentiable $\chi_i\colon\R\times (0,\infty)\to\R$.

The terms of the Jacobian matrix of  $\bW(\bc)$ can be computed as follows:
\begin{equation*}   \frac{\partial W_i(\bc)}{\partial c_j}=\delta_{i,j}\frac{\partial \chi_i}{\partial \xi}(c_i, \varphi(b^T\bc)) +
\frac{\partial \chi_i}{\partial p} (c_i, \varphi(b^T\bc)) b_{j}\varphi'(b^T\bc),
\end{equation*}
where $\delta_{i,j}$ is the Kronecker delta symbol.
Therefore, $\bW'(\bc)$ has the structure
\begin{align}\label{eq:WprimaD+AB}
\left\{\begin{aligned}
&\bW'(\bc) = D(\bc) + B(\bc) A(\bc)^T, \\
&D(\bc)=\text{diag}(v_i(\bc)), \quad v_i(\bc)={\frac{\partial \chi_i}{\partial \xi}(c_i, \varphi(b^T\bc))},\,\, i=1,\dots,N, \\
&B_i(\bc)=  \displaystyle{\frac{\partial \chi_i}{\partial p} (c_i, \varphi(b^T\bc))},
  \quad A_i(\bc)= b_{i}\varphi'(b^T\bc),\,\, i=1,\dots,N.
\end{aligned}\right.
\end{align}

The following result (see \cite{DGM18, Donat3, Anderson}) states that, under certain conditions, the eigenstructure of such matrices can be obtained from the roots of  the  function
\begin{equation*}
  Q(\lambda)=Q[\bc](\lambda)=1 + \sum_{j=1}^{N}
\frac{\gamma_j(\bc)}{v_j(\bc) -\lambda}, \quad \gamma_j(\bc)=A_j(\bc)B_j(\bc).
\end{equation*}

\begin{proposition}\label{pr:secular}
Let  $v,A,B\in\R^N$ and
\begin{align*}
    Q(\lambda)=1+\sum_{j=1}^{N}\frac{A_jB_j}{v_j-\lambda}.
\end{align*}
If $v_1  < v_2 < \cdots < v_N , \,\, A_i  B_i <0, \, 1 \leq i \leq  N$,
and there exists $\lambda_*<v_{1}$ such that $ Q(\lambda_*)> 0$, then $Q$ has exactly $N$ roots $\lambda_1,\dots,\lambda_N$, such that
\begin{equation*}\lambda_* < \lambda_1 < v_1  <\lambda_2 < v_2 < \cdots  <\lambda_N <  v_N ,
\end{equation*}
and these are all the eigenvalues of $M=D+BA^T$, $D=\text{diag}(v_1,\dots,v_N)$, which is therefore diagonalisable, and the (right)  eigenvectors can be explicitly determined from the eigenvalues as $r_j=(D-\lambda_j I_{N})^{-1}B$, i.e.,
\begin{equation}\label{eq:eigenvectors}
 (r_j)_k=\frac{B_k}{v_{k}-\lambda_j},\quad j,k=1,\dots,N.
\end{equation}
\end{proposition}

\begin{theorem} \label{th:Cw-welldefToth}
The function $\bW\colon (0,\infty)^{N}\to(0,\infty)^{N}$ in \eqref{eq:Wi-defToth} is a continuously differentiable diffeomorphism satisfying that $\bW'(\bc)$  is diagonalizable with  $N$ distinct eigenvalues $\lambda_i(\bc),\,\, i=1,\dots,N$,
\begin{equation}\label{eq:interlacing}
1 < \lambda_1(\bc) < v_1 (\bc) <\lambda_2(\bc) < v_2(\bc) < \cdots
<\lambda_N(\bc) <  v_N (\bc)\leq 1+\eta_{N},
\end{equation}
and eigenvectors explicitly computable from each $\lambda_i(\bc)$.

For $\bC=\bW^{-1}$, the eigenvalues of $\bC'(\bw), \bw\in(0,\infty)^N$, are
$$\mu_{j}(\bw)=\frac{1}{\lambda_j(\bC(\bw))},\quad j=1,\dots,N,$$
and satisfy
\begin{equation}\label{eq:interlacinv}
1 > \mu_1(\bw)> \mu_2(\bw)> \cdots  >\mu_N(\bw)> \frac{1}{1+\eta_N}.
\end{equation}
The corresponding eigenvectors are those of $\bW'(\bC(\bw))$, given by \eqref{eq:eigenvectors}.
\end{theorem}

\begin{proof}
The function $\bW(\bc)$ in \eqref{eq:Wi-defToth} is  given by:
\begin{align*}
  W_i(\bc)&=\chi_i(c_i, \varphi(b^T\bc)),\\
\chi_{i}(\xi, p)&=\xi(1+\frac{\eta_i}{p}),
\end{align*}
where $0<\eta_1<\eta_2<\dots<\eta_N$ and $b_i>0$. Therefore, for any $ \bc \in (0,\infty)^N$,  $\bW'(\bc)$ takes the form \eqref{eq:WprimaD+AB}, with
\begin{align*}
 v_i(\bc)&=\ffc{i}{\varphi(b^T\bc)} = 1+\frac{\eta_i}  {\varphi(b^T\bc)},  \\[5pt]
B_i&=\ffp{i}{\varphi(b^T\bc)}=-\frac{c_i\eta_i}{\varphi(b^T\bc)^2}, \\   A_i&=b_i\varphi'(b^T\bc).
\end{align*}
Since $0<\eta_1<\eta_2<\dots<\eta_N$ and $\varphi(b^T\bc)\geq 1$,
\begin{equation}\label{eq:interlacdir}
  1 <v_1(\bc) <v_2(\bc) < \cdots < v_N(\bc)\leq 1+\eta_{N}.
\end{equation}  
Furthermore, $\bc\in(0,\infty)^N $ yields that $b^T\bc > 0$ and $\varphi'(b^T\bc) > 0$  by hypothesis, so that $A_iB_i < 0,\forall i=1,\dots,N$.

We aim to apply Proposition \ref{pr:secular} with $\lambda_*=1$, for which we only have to check that $Q(1) > 0$:
\begin{align*}
  Q(1)&=1+\sum_{i=1}^{N} \frac{A_{i}B_{i}}{v_i(\bc)-1} =
    1-\sum_{i=1}^{N} \frac{c_i}{\varphi(b^T\bc)}b_i\varphi'(b^T\bc)\\
  &=\frac{\varphi(d)-d \varphi'(d)}{\varphi(d)}=\varphi(d)
    \left(\frac{d}{\varphi(d)}\right)' > 0, \,\, d=b^T\bc>0.
\end{align*}
Proposition \ref{pr:secular} thus yields that the eigenvalues of $\bW'(\bc)$ are $\lambda_1(\bc),\dots,\lambda_{N}(\bc)$ and satisfy the interlacing property \eqref{eq:interlacing}.

For any $\bw\in(0,\infty)^N$, let $\bc\in(0,\infty)^N$ be such that $\bW(\bc)=\bw$. Since the eigenvalues of $\bW'(\bc)$ are non-zero, the Inverse Function Theorem implies that $\bW$ is a local bijection with continuously differentiable local inverse. Since $\bW$ is already (globally) bijective, it follows that  $\bC=\bW^{-1}\colon (0,\infty)^N\to(0,\infty)^N$  is continuously differentiable.

The chain rule implies $\bC'(\bw)=(\bW'(\bC(\bw)))^{-1}$, hence, the eigenvalues of $\bC'(\bw)$ are  $1/\lambda_i(\bC(\bw))$, $i=1,\dots,N$, and its eigenvectors are those of $\bW'(\bC(\bw))$. Therefore,  \eqref{eq:interlacinv} is a direct consequence of \eqref{eq:interlacing}.
\end{proof}

From Theorem \ref{th:Cw-welldefToth},  system \eqref{eq:modelwc1} can be rewritten as:
\begin{equation} \label{eq:modelwf}
\frac{\partial {\boldsymbol w}}{\partial t}+\frac{\partial (u
  {\boldsymbol C}(\boldsymbol w))}{\partial z}=
\frac{\partial}{\partial z}\left[D_a\boldsymbol C'(\boldsymbol
  w)\frac{\partial {\boldsymbol w}}{\partial z}\right].
\end{equation}

\begin{corollary}\label{cor:1}
The system of conservation laws \eqref{eq:modelwf} is strictly hyperbolic in $\Omega=(0,\infty)^{N}$ when $D_a=0$. Moreover, the eigenvalues of the Jacobian matrix $\boldsymbol{f}'(\bw)$ are non-negative, pairwise distinct, and bounded above by $u$, for any $\bw\in\Omega$. If $D_a>0$, system \eqref{eq:modelwf} is parabolic in the sense of Petrovskii (cf. \cite{EidelmanZhiharatshu98}), i.e., there exists a positive lower bound of the eigenvalues of the diffusion matrix $D_a\bC'(\bw)$, for all $\bw\in\Omega$.
\end{corollary}

\begin{proof}
    Theorem \ref{th:Cw-welldefToth} ensures that  the eigenvalues $\mu_{j}(\bw)$ of  $\boldsymbol{C}'(\boldsymbol{w})$  satisfy \eqref{eq:interlacinv}.
    
Since  the convective flux is given by $ \boldsymbol{f(w)} = u \boldsymbol{C}(\boldsymbol{w})$, its Jacobian matrix is  $\boldsymbol{f'(w)} = u \boldsymbol{C}'(\boldsymbol{w})$, and its eigenvalues,  $u\mu_j(\bw)$, $j=1,\ldots, N$, are positive, pairwise distinct and bounded above by $u$.

For $D_a>0$, the eigenvalues of the diffusion  matrix $D_a\boldsymbol C'(\boldsymbol w)$ are
$D_a\mu_j(\bw)$, $j=1,\ldots, N$, and are bounded below by $\frac{D_a}{1+\eta_N} > 0$.
\end{proof}

\section{Characteristic-based WENO schemes}
\label{sec:charweno}

Essentially non-oscillatory high-order schemes for hyperbolic systems of conservation laws need to use reconstructions whose stencils avoid discontinuities as much as possible. This has been satisfactorily tackled for scalar conservation laws by several methods, such as the MUSCL schemes from the pioneering work of van Leer \cite{vanLeer79}, ENO schemes \cite{EngquistHartenOsher86} and WENO schemes \cite{LiuOsherChan94, JiangShu1996, Borges2008, Castro2011,Rat1, Rat2, Rat3, Rat4, Rat5}. Since one of the purposes of this paper is to show that local characteristic projections can be used in this context to prevent non-physical oscillations, we use standard WENO techniques, with Jiang-Shu's weights \cite{JiangShu1996},  although more advanced WENO techniques, as WENO-Z \cite{Borges2008}, could be used as well. The common feature of these reconstructions is that they achieve their accuracy while reducing  numerical oscillatory behaviour if discontinuities are well separated, a fact that will not generally occur for systems,  where waves associated with different characteristic fields may interact.

For linear hyperbolic systems, one can resort to changing to characteristic variables, by using projections onto (left) eigenvectors of the flux matrix, thus getting a fully decoupled system of (scalar) conservation laws which would be amenable to oscillations-free reconstructions.

In the nonlinear case, there is no change of variables that would decouple the system, let alone maintain weak solutions. What is heuristically used are {\em local characteristic variables and fluxes}, which consist in computing numerical fluxes at some cell interface by using reconstructions of projections onto (left) eigenvectors of a Jacobian matrix of the flux associated with that cell interface.

It is well known that centred schemes are unstable for hyperbolic equations, so another difficulty faced in the design of numerical schemes for hyperbolic systems of conservation laws is that of maintaining stability by appropriate {\em upwinding}, which amounts to explicitly or implicitly adding  {\em numerical viscosity}.

Since Corollary \ref{cor:1} yields that  characteristic velocities, i.e., the eigenvalues of  Jacobian matrices of the flux, are positive, the natural choice is to use left-biased reconstructions, which are functions $\mathcal{I}$ of an odd number $2s+1$ of arguments, such that
\begin{align*}
  \mathcal{I}(\overline{f}_{-s,h},\dots,\overline{f}_{s,h})=f(h/2) +  \bigO(h^p), \quad
  \overline{f}_{j,h}=\frac{1}{h}\int_{(j-\frac{1}{2})h}^{(j+\frac{1}{2})h} f(x)dx,
\end{align*}
for any  real function $f$ which is sufficiently smooth in a neighborhood of $0$ and where $p$ is the {\em order} of the reconstruction.

 We follow Shu and Osher's methodology for the computation of the numerical flux $\widehat{\boldsymbol{f}}_{j+\mig}$ and use the matrix $R$ whose columns are the eigenvectors $R_k=r_{k}(\mig(\bw_{j}+\bw_{j+1}))$ obtained in \eqref{eq:eigenvectors} from  Theorem \ref{th:Cw-welldefToth} and Proposition \ref{pr:secular} to compute local characteristic fluxes
\[\widetilde{\boldsymbol{f}}_{k,l}=(R^{-1}\boldsymbol{f}(\bw_{k}))_{l},\]
with $k=j-s,\dots,j+s$, $l=1,\dots,N$, and
\begin{align}\label{eq:charflux}
  \widehat{\boldsymbol{f}}_{j+\mig}=
  \sum_{l=1}^{N} R_{l} \mathcal{I}(\widetilde{\boldsymbol{f}}_{j-s,l},\dots,\widetilde{\boldsymbol{f}}_{j+s,l}).
\end{align}

\section{IMEX schemes}
\label{sec:imex-weno}

In this section, we review the numerical scheme that was proposed in
\cite{DGM18} (to which we refer  the reader for specific details) for the solution of  the conservative formulation of the
ED model \eqref{eq:q-def}-\eqref{eq:modelwc1},  which can be  rewritten as follows:
\begin{equation}
\label{eq:modelwc11}
\left\{\begin{aligned}
&\frac{\partial \bw}{\partial t} + \frac{\partial}{\partial z}\left(\boldsymbol{f}(\bw)- \boldsymbol{g}\left(\bw, \frac{\partial \bw}{\partial z}\right)\right)=0, \\
&\boldsymbol{g}(\bw, \frac{\partial \bw}{\partial z})=D_a \frac{\partial  \bC(\bw)}{\partial z}=D_a \bC'(\bw)\frac{\partial \bw}{\partial z},
\end{aligned}
\right.
\end{equation}
provided with corresponding left (up) Dankwerts and right (down) Neumann boundary conditions
\begin{equation} \label{eq:bc2w}
(\boldsymbol{f}(\bw)-\boldsymbol{g}\left(\bw, \frac{\partial \bw}{\partial z}\right))(0,t)=u \bc_{inj}(t), \qquad  \boldsymbol{g}\left(\bw, \frac{\partial \bw}{\partial z}\right)(1,t)= 0,
\end{equation}
where $\bc_{inj}(t)$ is the vector of concentrations of injected
components at time $t$.

The scheme is obtained by the Method Of Lines and yields, in the first stage,  approximations which are the solutions of an ODE system
\begin{equation} \label{eq:weno-sd}
\boldsymbol{w}'(t)=\boldsymbol{{\cal{L}}}(\bw(t), t)+ \boldsymbol{{\mD}}(\bw(t)),
\end{equation}
where $\bw(t)$ is an $N\times m$ matrix whose $j$-th column, $\bw_j(t)$, is an approximation of $ w(z_j, t)\in\mathbb R^N$, $z_j=(j-\mig)\Delta z, j=1,\dots, m, \Delta z=\frac{1}{m}$, and  $\boldsymbol{{\mD}}$ is the spatial discretization of the diffusion term $\frac{\partial }{\partial z}\boldsymbol{g}\left(\bw, \frac{\partial \bw}{\partial z}\right)$ in \eqref{eq:modelwc11} given by $\mD(\bw)=D_a\boldsymbol{C}^*(\bw)\mathcal{A}$,
where $\mathcal{A}$ is  the tridiagonal $m\times m$ matrix that
discretizes the 1D Laplacian operator with Neumann boundary conditions and we use the notation
\begin{align*}
\boldsymbol{W}^*(\bc)_{i,j}=W_i(\bc_j),\quad
\boldsymbol{C}^*(\bw)_{i,j}=C_i(\bw_j).
\end{align*}
Note that $\bw, \bc$ are $N\times m$ matrices and $\bw_j, \bc_j$, $j=1,\dots,m$,  are their respective $j$-th columns.

Finally, the  term $\boldsymbol{{\cal{L}}}$ is an $N\times m$ matrix whose $j$-th column
\begin{align}\label{eq:conv-diff}
\boldsymbol{{\cal L}}_j = -\frac{1}{\Delta z} \left(
  \boldsymbol{\widehat{{f}}}_{j+1/2}-\boldsymbol{\widehat{{f}}}_{j-1/2}\right),\qquad
\end{align}
is the spatial discretization of the convective term $-\frac{\partial}{\partial z} \boldsymbol{f}(\bw)(z_{j},t)$, obtained through the following numerical fluxes $\boldsymbol{\widehat{{f}}}_{j+1/2}$ (the explicit dependence of $\boldsymbol{{\cal{L}}}$ on $t$ is due to the boundary term containing $\bc_{\text{inj}}(t)$ in \eqref{eq:bc2w}):
\begin{enumerate}
\item  The characteristic-based numerical fluxes \eqref{eq:charflux},
  introduced in Section \ref{sec:charweno} for the fifth order WENO
  (WENO5) reconstruction $\mathcal{I}$ ($s=2$) with Jiang-Shu's weights  \cite{JiangShu1996}.

\item The first order upwind numerical fluxes
\begin{equation} \label{eq:compupw1}
\widehat{\boldsymbol{f}}_{j+\mig}=\boldsymbol{f}(\bw_{j}).
\end{equation}

\item The fifth order upwind WENO5 numerical fluxes
\begin{equation}\label{eq:compupw5}
(\widehat{\boldsymbol{f}}_{j+\mig})_l=\mathcal{I}\big(\boldsymbol{f}_{l}(\bw_{j-2}),\dots,\boldsymbol{f}_{l}(\bw_{j+2})\big),\,\, l=1,\dots,N.
\end{equation}

\item The fifth order  Global Lax-Friedrichs WENO5 numerical fluxes
\begin{align}\label{eq:compglf}
(\widehat{\boldsymbol{f}}_{j+\mig})_l&=\mathcal{I}(\boldsymbol{f}^+_{j-2,l},\dots,\boldsymbol{f}^+_{j+2,l})+\mathcal{I}(\boldsymbol{f}^-_{j+3,l},\dots,\boldsymbol{f}^-_{j-1,l}),\\
\notag
\boldsymbol{f}^{\pm}_{k,l}&=\frac{1}{2}(\boldsymbol{f}(\bw_{k})\pm \alpha \bw_{k})_{l},\,\, l=1,\dots,N,
\end{align}
where the numerical viscosity $\alpha$ is a local-in-time upper bound of all the characteristic speeds of the problem. For the ED model, Corollary \ref{cor:1} ensures that we can  take $\alpha=u$.
\item  Characteristic-based WENO5 numerical fluxes  with
  the Global Lax-Friedrichs flux-splitting
\begin{align}\label{eq:charglf}
  &\widehat{\boldsymbol{f}}_{j+\mig}=
  \sum_{l=1}^{N} R_{k} \Big(
\mathcal{I}(\widetilde{\boldsymbol{f}}^{+}_{j-2,l},\dots,\widetilde{\boldsymbol{f}}^{+}_{j+2,l})+ \mathcal{I}(\widetilde{\boldsymbol{f}}^{-}_{j+3,l},\dots,\widetilde{\boldsymbol{f}}^{-}_{j-1,l})\Big),\\
  \notag
&\widetilde{\boldsymbol{f}}^{\pm}_{k,l}=\frac{1}{2}\big(R^{-1}(\boldsymbol{f}(\bw_{k})\pm \alpha \bw_{k})\big)_{l},\,\, l=1,\dots,N,
\end{align}
where the $N\times N$ matrix $R$ is defined in Section \ref{sec:charweno} and $\alpha$ is as in the previous item.

\item The second-order upwind MUSCL numerical fluxes given by
\begin{align}\label{eq:musclupw}
  &\widehat{\boldsymbol{f}}_{j+\mig}=\boldsymbol{f}(\boldsymbol{w}_{j+\mig}),\\
  \notag
  &\boldsymbol{w}_{j+\mig}=\bw_{j}+\mig\text{minmod}(\bw_{j}-\bw_{j-1},
    \bw_{j+1}-\bw_{j}),\\
  \notag
  &\text{minmod}(a,b)=\frac{\text{sign}(a)+\text{sign}(b)}{2}\min(|a|,|b|).
\end{align}
\end{enumerate}

To obtain fully-discrete, high-order, schemes, an appropriate  ODE solver must be applied to \eqref{eq:weno-sd}. As mentioned in \cite{DGM18}, explicit Runge-Kutta solvers require that the time step $\Delta t$ be proportional to $\Delta z^2$ for stability, due to the explicit treatment of the second-order term. The benefits of using an implicit-explicit (IMEX) scheme, that treats implicitly only the term requiring the mentioned time step restriction, have already been assessed in \cite{DGM18} in the context of the conservative ED model with multicomponent Langmuir isotherms. We consider here the implicit-explicit midpoint rule (which is second-order in time, see, e.g., \cite{ARS97})
\begin{align}\label{imexrk21}
\boldsymbol{w}^{n+1/2}&=\boldsymbol{w}^n+\frac{\Delta t}{2}
\left(\boldsymbol{{\cal L}}(\bw^n, t_n)+\boldsymbol{{\mD}}(\bw^{n+1/2}) \right), \\
\boldsymbol{w}^{n+1}&=\boldsymbol{{w}}^n+\Delta t
\left(\boldsymbol{{\cal L}}(\bw^{n+1/2}, t_{n+\mig})+\boldsymbol{{\mD}}(\bw^{n+1/2})
\right), \label{imexrk22}
\end{align}
where  $t_{n+1/2}=t_n+\Delta t_n/2$.

Note that the second step \eqref{imexrk22} is explicit, thus ${\boldsymbol w}^{n+1}$ can be directly obtained from $\boldsymbol{w}^{n+1/2}$ and $\boldsymbol{w}^{n}$.
To solve \eqref{imexrk21},
we can express the terms depending on $\bw^{n+1/2}$ in terms of $\bc^{n+1/2}=\bC^*(\bw^{n+1/2})$ (taking into account that $\bw^{n+1/2}=\bW^*(\bc^{n+1/2})$) to arrive at the matrix equation:
\begin{align*}
&\boldsymbol{W}^*(\bc^{n+1/2})-\frac{D_a\Delta t}{2} \bc^{n+1/2}\mathcal{A} =\mathcal{G}^n,\quad
\mathcal{G}^n= \bw^n+\frac{\Delta t}{2} \boldsymbol{{\cal L}}(\bw^n),
\end{align*}
which, dropping the $n+\mig$ superindex, is equivalent to
\begin{align}\label{eq:pep1}
&c_{i,j}\left(1+\frac{\eta_{i}}{\varphi(b_1c_{1,j}+\dots+b_Nc_{N,j})}\right)-\frac{D_a\Delta t}{2}
(\bc\mathcal{A})_{i,j}=\mathcal{G}_{i,j}^n,\end{align}
where  $i\in\{1, \dots , N\}$ refers to the component of the mixture and $j\in\{1, \dots , m\}$ refers to the grid point under consideration.

Therefore, solving \eqref{eq:pep1} by Newton's method involves solving a block-tri\-diago\-nal system with small $N \times N$ blocks at each iteration step.   A standard block tridiagonal LU factorization algorithm (see \cite{Golub} for details) can be used to solve this efficiently.

\section{Numerical experiments} \label{sec:numex}
In this section, we perform some numerical experiments to compare the performance of the proposed IMEX scheme coupled with the reconstruction procedure and numerical fluxes introduced in Sections \ref{sec:charweno} and \ref{sec:imex-weno}, for the ED model with T\'oth's isotherms.

In the remainder of the section we will use the following notation: by CHR-UPW we will denote the IMEX scheme coupled with the characteristic-based numerical fluxes defined by \eqref{eq:charflux}. We use COMP-UPW1 to refer to the IMEX scheme that uses the first order upwind numerical fluxes in \eqref{eq:compupw1} and, analogously, COMP-UPW5 will refer to the same scheme, but using the fifth order upwind WENO5 numerical fluxes \eqref{eq:compupw5}. COMP-GLF will be used for the IMEX scheme coupled with the Global Lax-Friedrichs numerical fluxes in \eqref{eq:compglf}, proposed in \cite{DGM18}. The Global Lax-Friedrichs numerical fluxes can also be coupled with the characteristic information of the system, as explained in Section \ref{sec:imex-weno}. We use the notation CHR-GLF for the IMEX scheme coupled with the numerical fluxes given by \eqref{eq:charglf}. Finally, we denote the IMEX scheme with the second-order upwind MUSCL numerical fluxes given in \eqref{eq:musclupw} by MUSCL.

To compare the numerical results obtained using different schemes, we will compute the experimental approximate $L^1$-errors. Denoting by $\smash{(w_{i,j}^m(t))_{j=1}^{m}}$ and $\smash{(w_{i,l}^{\mathrm{ref}}(t))_{l=1}^{m_{\mathrm{ref}}}}$ the numerical solution for the $i$-th component at time $t$ calculated with $m$ and $m_{\mathrm{ref}}$ cells, respectively, we compute $\smash{\tilde{w}_{i,j}^{m}(t)}$ for $j=1,\dots,m$ by
\[ \tilde{w}^{m}_{i,j}(t)=\frac{1}{ R}\sum_{k=1}^{ R} w^{\mathrm{ref}}_{i,R(j-1)+k}(t), \quad R=m_{\mathrm{ref}}/m.\]
Then, the total approximate $L^1$-error of the numerical solution  $\smash{(w_{i,j}^m(t))_{j=1}^{m}}$ at  time~$t$ is  then given  by
\begin{align} \label{formula_error}
e_m (t)=\frac{1}{m} \sum_{i=1}^N \sum_{j=1}^{m} \bigl|\tilde{w}_{i,j}^{m}(t)-w_{i,j}^{m}(t)\bigr|.
\end{align}

\subsection{Experiment 1}

In this experiment, we model the elution chromatography process considering three components proposed in \cite{Javeed} (Section 4.3). As we have previously mentioned, elution chromatography is a separation process used to separate a mixture of compounds using a solid stationary phase and a liquid moving phase. The process takes place in a column filled with the stationary phase where we introduce the mixture that we want to separate along with a mobile phase that moves the sample mixture through the column.

We consider here a mixture sample that contains two solutes, components 1 and 2. A third component, known as the mobile phase, the solvent or the displacer, is also pumped into the column at a certain flow rate. The descent of components 1 and 2 through the column is determined by their relative affinities for the mobile and stationary phases. The main characteristic of this process is that the components of the mixture can separate completely inside the column, forming bands of high concentration of one component as the mixture goes down the column. The series of such bands is called the \textit{isotachic train} (\cite{cazes2001encyclopedia}).

The parameters that we consider for this experiment are $a_1=4, a_2=5, a_3=6, b_1=4, b_2=5, b_3=1$. Components 1 and 2 are only initially injected between $t=0$ and $t=0.1$ with $c_1=c_2=1\,\mathrm{g/l}$ at the top of the column, i.e. at $z=0\,\mathrm{m}$, while component 3, the displacer, is injected continuously after the mixture, from $t=0.1$, with $c_3=1\,\mathrm{g/l}$. In addition, we consider $\epsilon=0.5$ and $u=0.2$.

As stated in \cite{DGM18}, the stability restrictions for the proposed IMEX scheme \eqref{imexrk21}-\eqref{imexrk22} are given by
\begin{equation} \label{IMEXstability}
\frac{u\Delta t}{\Delta z}\max_{w} \varrho (\bC'(w))=K\leq 1,
\end{equation}
where $\varrho (\bC'(w))$ denotes the spectral radius of the Jacobian matrix $\bC'(w)$. For the characteristic-based schemes, we use the computed value of the eigenvalues to determine the maximum in \eqref{IMEXstability} and obtain the value of the time step $\Delta t$ in each iteration, for a fixed value of $\Delta z$. For component-wise schemes, the spectral information of the Jacobian matrix is not computed. Therefore, using that $\max_{w} \varrho (\bC'(w))\leq 1$, we will use the condition
\begin{equation*} u \frac{\Delta t}{\Delta z}=K\leq 1,
\end{equation*}
to determine the value of $\Delta t$. For all the numerical experiments in this work, we will use  $K=0.8$.

\begin{figure}
\includegraphics[width=\linewidth]{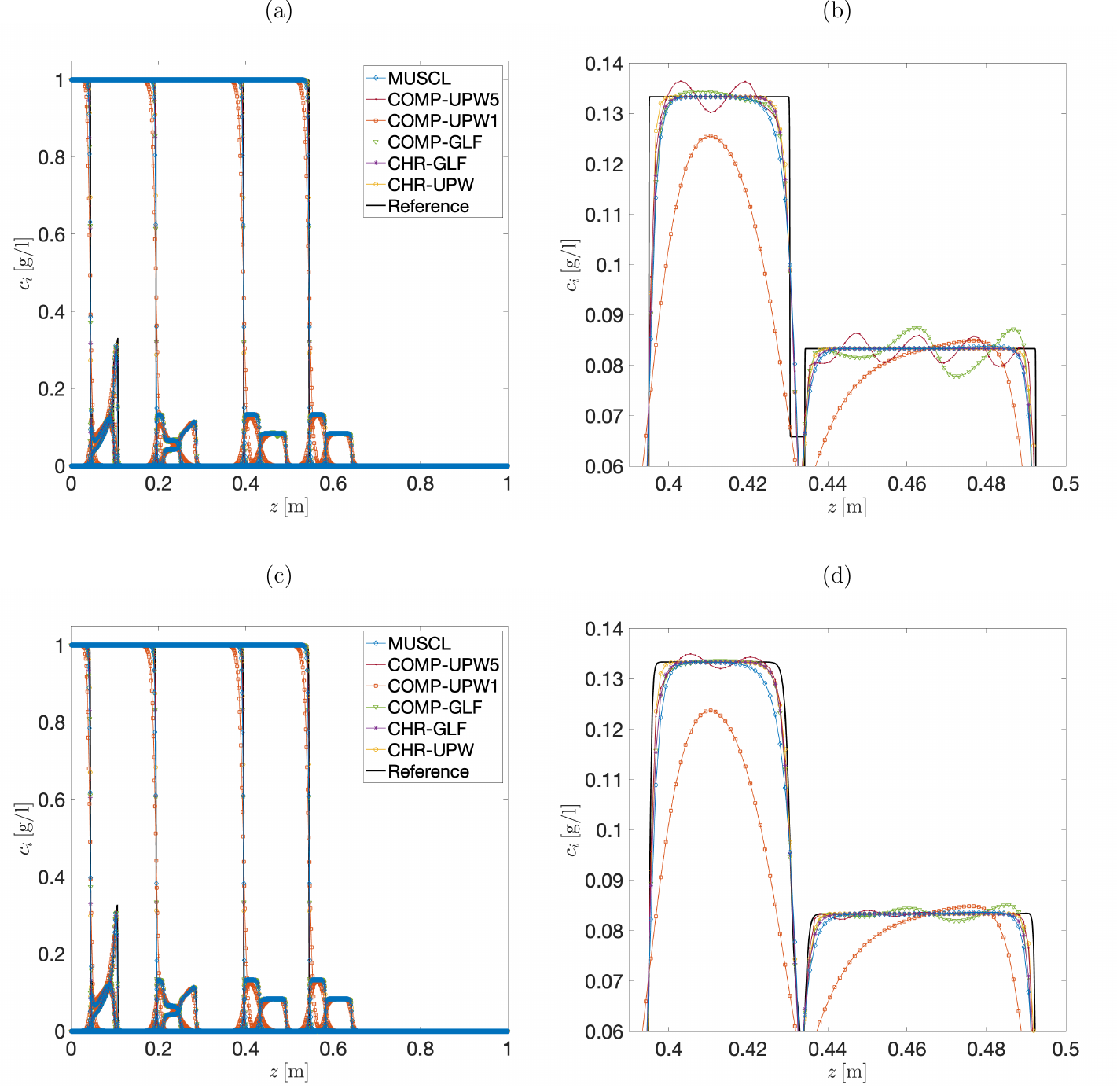}
\caption{Experiment 1. Numerical solutions obtained with MUSCL, COMP-UPW5, COMP-UPW1, COMP-GLF, CHR-GLF and CHR-UPW schemes for $\nu=1$, $D_a=0$ (a) and $D_a=10^{-5}$ (c) for  $T=1, 4, 8$ and $11$. Enlarged views of approximate concentrations  for components 1 and 2 at $T=8$ are given in (b) and (d).}
\label{comp3disperso1}
\end{figure}

In Figure \ref{comp3disperso1} (a) and (c), we show the approximate solutions obtained for $D_a=0$ and $D_a=10^{-5}$, respectively, a mesh with $m=800$ and $\nu=1$, i.e., we use the {\it Langmuir} isotherms. We show the approximate solutions for $T=1, 4, 8$ and $11$, where the formation of the isotachic train can be clearly appreciated. The reference solution is also included in the plots. It is computed with the CHR-UPW scheme and a computational mesh of $m=25600$ nodes.

Figure \ref{comp3disperso1} (b) and (d) correspond to enlarged views of the recovered solutions for $T=8$ and the two values of $D_a$ chosen. As can be seen, when a steep step structure appears in the solution, the numerical solutions obtained using COMP-GLF and COMP-UPW5 schemes show spurious oscillations that are not present when using schemes that use the characteristic information of the system. Approximate solutions obtained with COMP-GLF and COMP-UPW5 have similar behaviour (in terms of the amplitude of the oscillations and accuracy near the discontinuity), while the solution obtained with the COMP-UPW1 method is less accurate, producing smoothed-out profiles. Finally, although the numerical solutions obtained using the MUSCL scheme do not present spurious oscillations, their accuracy, especially near sharp edges, is worse than the accuracy of characteristic-based schemes. For $D_a=10^{-5}$, even though the profiles of the numerical solutions are smoothed out, see Figure \ref{comp3disperso1} (d), some oscillatory behaviour can still be seen when using COMP-GLF and COMP-UPW5. It is also appreciated that CHR-GLF is more diffusive than CHR-UPW.

\begin{figure}
\includegraphics[width=\linewidth]{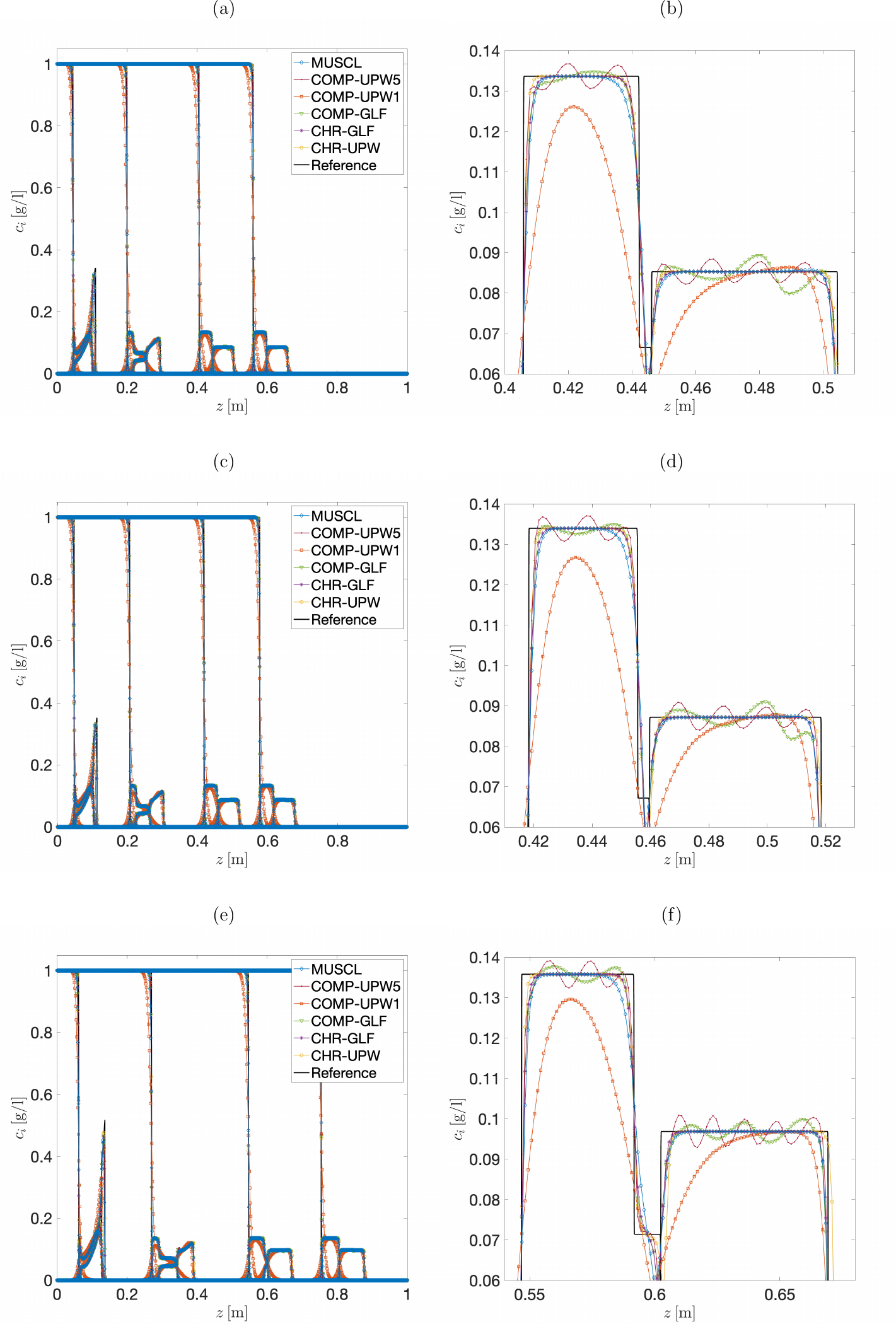}
\caption{Experiment 1. Numerical solutions obtained with MUSCL, COMP-UPW5, COMP-UPW1, COMP-GLF, CHR-GLF and CHR-UPW schemes with $D_a=0$ and $\nu=0.95$ (a), $\nu=0.9$ (c) and $\nu=0.6$ (e) for  $T=1, 4, 8$ and $11$. Plots (b), (d) and (f) are enlarged views of (a), (c) and (e) respectively, for  $T=8$.}
\label{comp3disperso2}
\end{figure}

\begin{figure}
\includegraphics[width=\linewidth]{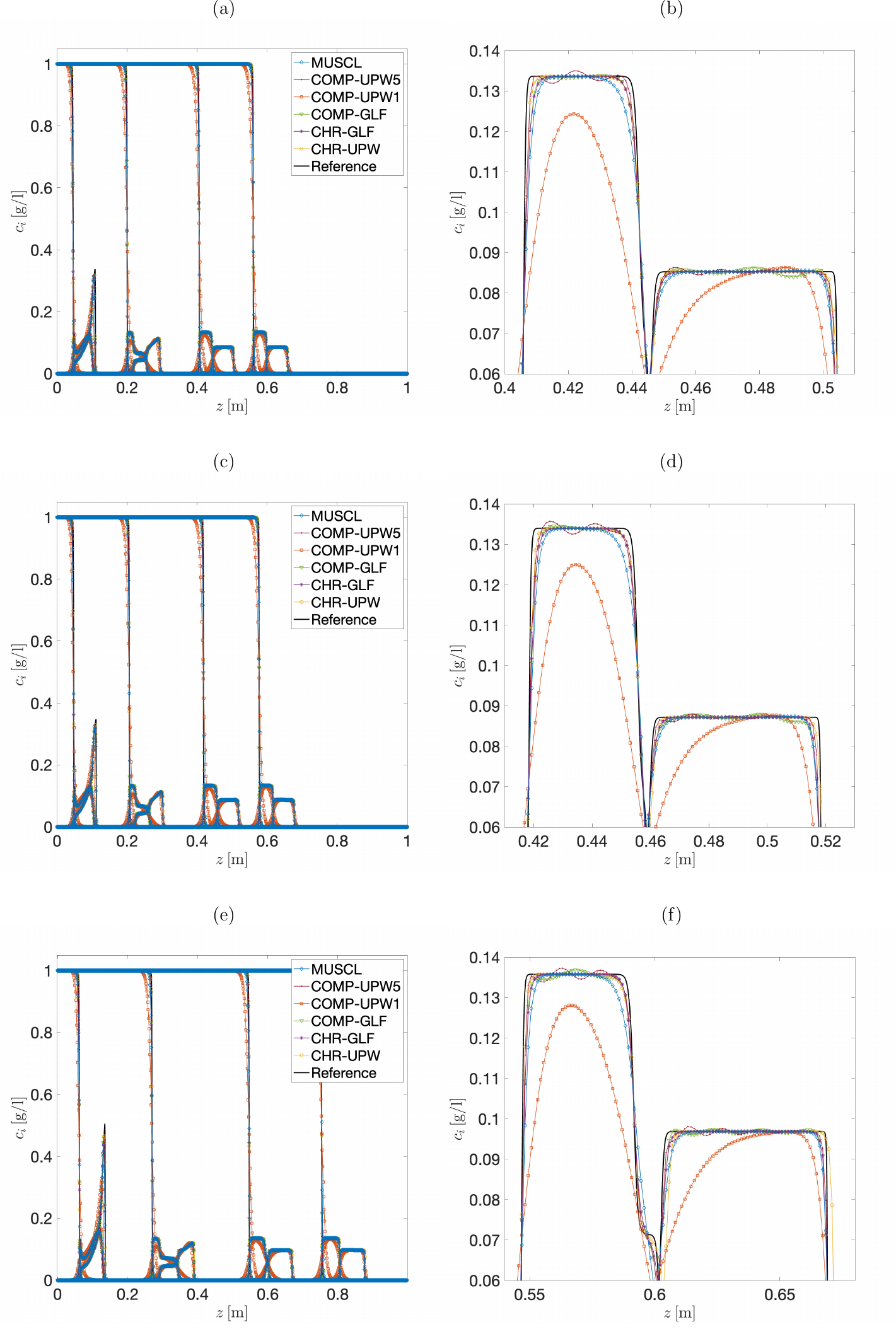}
\caption{Experiment 1. Numerical solutions obtained with MUSCL, COMP-UPW5, COMP-UPW1, COMP-GLF, CHR-GLF and CHR-UPW schemes with $D_a=10^{-5}$ and $\nu=0.95$ (a), $\nu=0.9$ (c) and $\nu=0.6$ (e) for  $T=1, 4, 8$ and $11$. Plots (b), (d) and (f) are enlarged views of (a), (c) and (e) respectively, for  $T=8$.}
\label{comp3disperso3}
\end{figure}

In Figures \ref{comp3disperso2} and \ref{comp3disperso3} the approximate solutions of the same scenario but with values of $\nu<1$ are displayed. In particular, we show the results for $\nu=0.95$ (a), $\nu=0.9$ (c) and $\nu=0.6$ (e) for $D_a=0$ in Fig. \ref{comp3disperso2} and $D_a=10^{-5}$ in Fig. \ref{comp3disperso3}. It is clearly appreciated that the smaller the value of $\nu$ considered, the faster the three components move through the column, which is in accordance with the  assumption that smaller values of $\nu$ correspond to a more heterogeneous stationary phase. The enlarged views of the approximate solutions for $T=8$ in both Figs. \ref{comp3disperso2} and \ref{comp3disperso3} (b), (d) and (f) show that the oscillatory behaviour observed in Figure \ref{comp3disperso1} for $\nu=1$ is also a feature that appears when we consider values of $\nu<1$, for both values of the parameter $D_a$ considered.

\begin{figure}
\includegraphics[width=\textwidth]{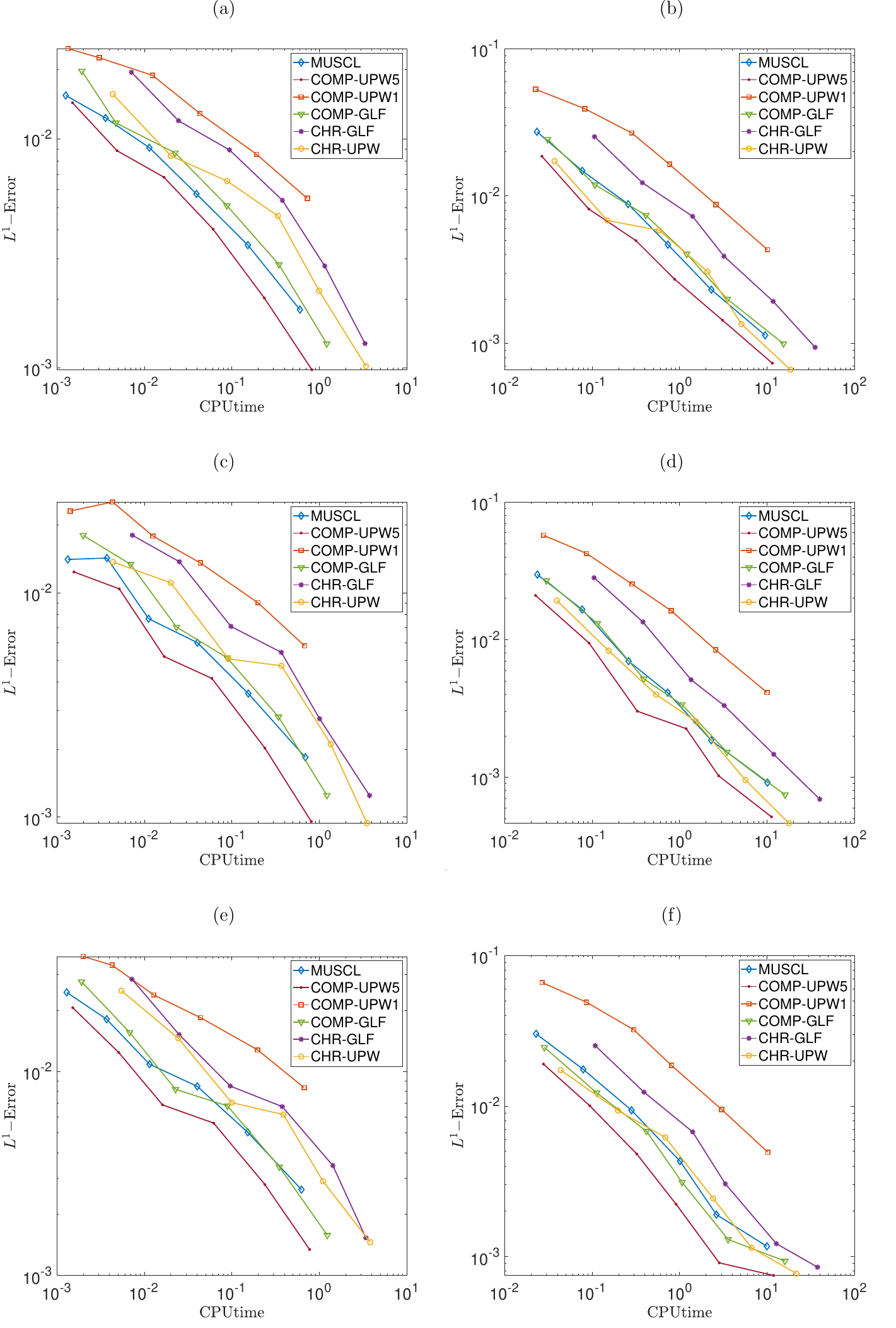}
\caption{Experiment 1. Performance of MUSCL, COMP-UPW5, COMP-UPW1, COMP-GLF, CHR-GLF and CHR-UPW methods for $D_a=0$, $T=1$ (left) and $T=11$ (right) and values of $\nu=1$ (a) and (b), $\nu=0.9$ (c) and (d) and $\nu=0.6$ (e) and (f).}
\label{comp3disperso4}
\end{figure}

To study the performance of each numerical method, we have represented in Figure \ref{comp3disperso4} the approximate $L^1$-errors vs. the computational times needed to obtain the approximate solutions with all the methods considered in this work, using a logarithmic scale. We have used $D_a=0$, three values of the parameter $\nu$ ($\nu=0.6,\, 0.9,\, 1$) and times $T=1$ and $11$.

As  can be seen, for $T=1$, the COMP-UPW5 scheme is the most efficient, since it obtains the most accurate approximate solutions in the shortest amount of computational time, for each value of the parameter $\nu$ considered. The accuracy of the CHR-UPW scheme is comparable but the computational time needed to obtain an approximate solution using this scheme is higher, since it is penalized by the numerical solution of the eigenvalue/eigenvector problems.

However, when a longer time is considered, for instance, $T=11$, the approximate solutions obtained using the CHR-UPW method show a better performance, closer to the one of the COMP-UPW5 scheme, as can be seen in Figure \ref{comp3disperso4} (b), (d) and (f). This feature can be explained by the fact that, for $T=11$, the isotachic train has fully developed and the oscillations in the approximate solutions obtained with component-wise methods decrease their accuracy.

\begin{figure}
\includegraphics[width=\textwidth]{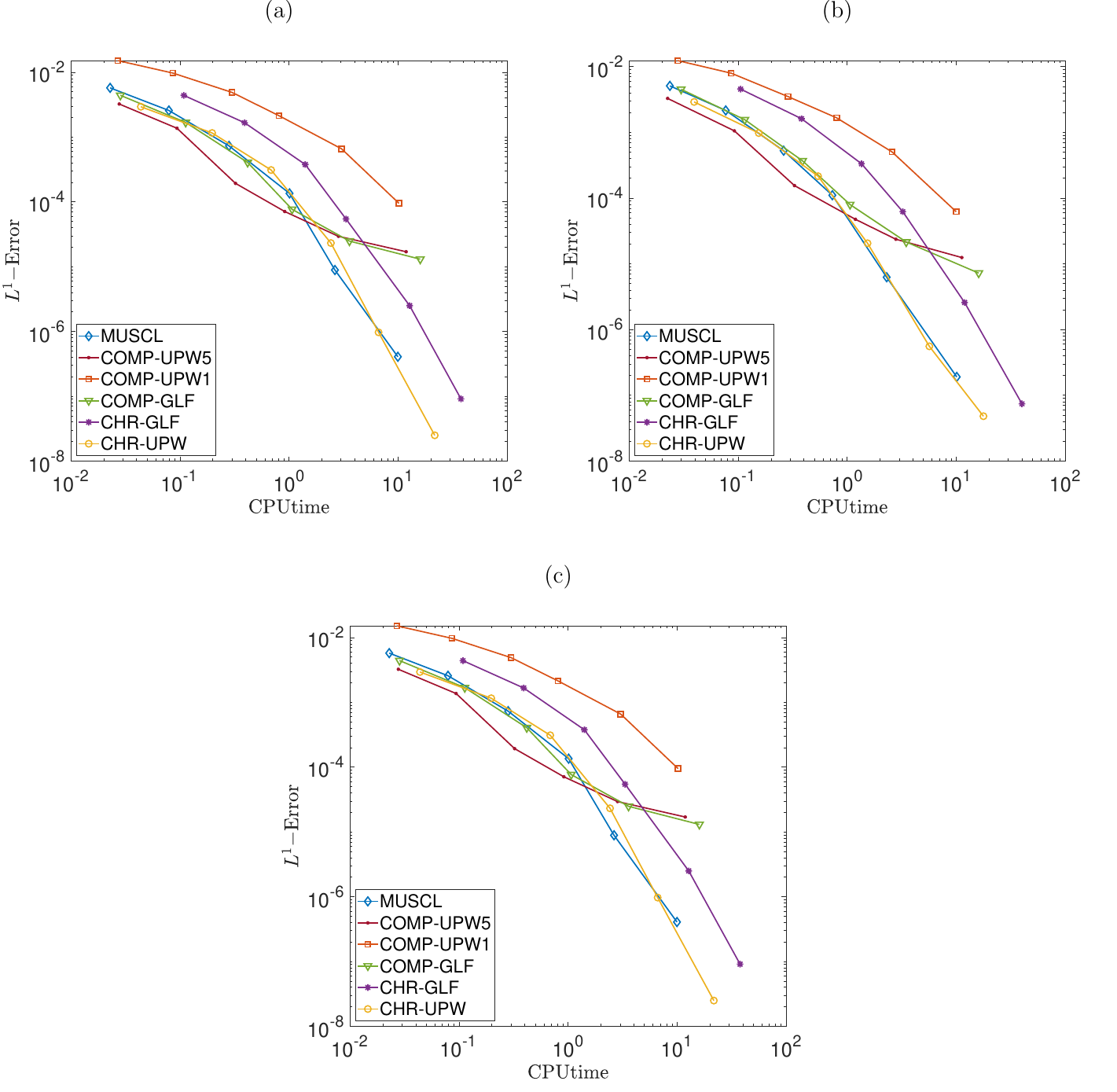}
\caption{Experiment 1. Performance of MUSCL, COMP-UPW5, COMP-UPW1, COMP-GLF, CHR-GLF and CHR-UPW methods, with the $L^1$-error computed discarding the $2\%$ of the largest errors. We have used $D_a=0$, $T=11$ and values of $\nu=1$ (a), $\nu=0.9$ (b) and $\nu=0.6$ (c).}
\label{comp3disperso5}
\end{figure}

From Figures \ref{comp3disperso1} - \ref{comp3disperso3}, it is clear that the approximate $L^1$-errors are dominated by the errors produced at the shocks. For $T=11$, we can repeat the efficiency plots shown in Figure \ref{comp3disperso4}, but leave out the largest errors, treating them as outliers, to better appreciate the performance of the CHR-UPW scheme in the oscillatory regions of the numerical solutions. In Figure \ref{comp3disperso5} we have discarded the $2\%$ of the largest errors and computed the approximate $L^1$-errors with the remaining $98\%$ of the values. We have then produced the plots in Figure \ref{comp3disperso4}comparing the approximate $L^1$-errors obtained with the computational times needed to obtain the approximate solutions with all schemes considered in this work and values of $\nu=0.6,\, 0.9,\, 1$.

As can be seen, the efficiency of both CHR-UPW and CHR-GLF schemes increases as the computational mesh size does. Even though the efficiency of MUSCL improves, the values of the error obtained are far from the ones obtained when using characteristic-based schemes. It is clear that the errors far away from shocks decrease faster for the schemes CHR-UPW and CHR-GLF than for the schemes  COMP-UPW  and COMP-GLF, which is in accordance with the results shown in Figures \ref{comp3disperso1} - \ref{comp3disperso3}.

\subsection{Experiment 2}

For this experiment, we consider the same parameters as in the previous one, but we diminish the quantity of displacer injected into the column. In particular, we consider $c_3=0.5 \,\mathrm{g/l}$.

\begin{figure}
\includegraphics[width=\linewidth]{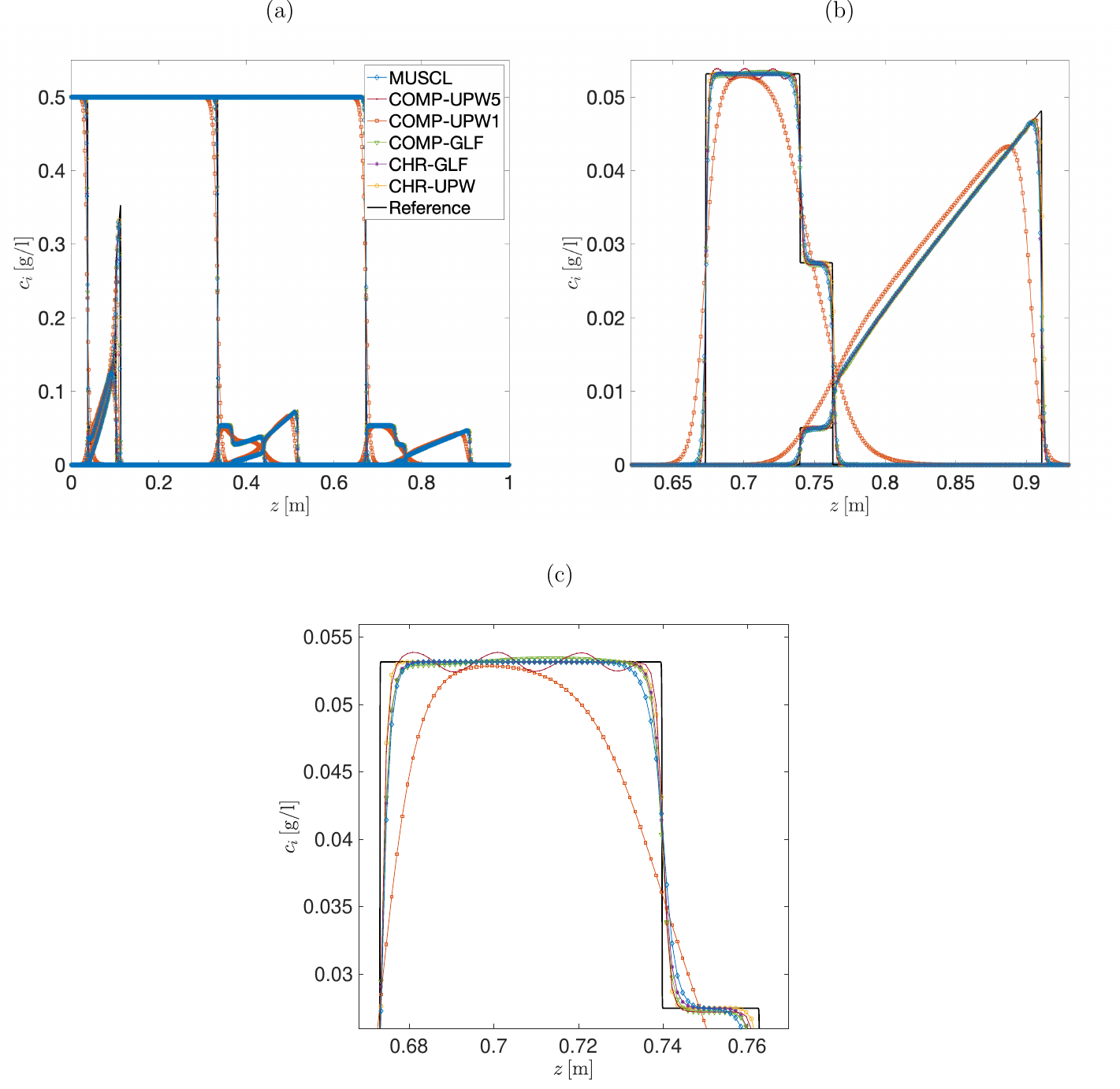}
\caption{Experiment 2. (a) Numerical solutions obtained with MUSCL, COMP-UPW5, COMP-UPW1, COMP-GLF, CHR-GLF and CHR-UPW schemes with $D_a=0$ and $\nu=0.9$ for  $T=1, 8$ and $16$.  (b) Enlarged view of the numerical solution of components 1 and 2 in (a) for  $T=16$.}
\label{comp3disperso6}
\end{figure}

In Figure \ref{comp3disperso6} we show the approximate solutions obtained for $T=1, 8$ and $16$ and $\nu=0.9$. As can be seen, components 1 and 2 do not separate completely in this scenario, as happened in the previous experiment, even if we run the experiment for a longer time. The enlarged views of the approximate solutions for $T=16$ in Fig. \ref{comp3disperso6} (b) show that, despite this fact, some numerical oscillations still appear in the numerical solutions near strong gradients when using component-wise schemes. These oscillations can be clearly appreciated in the numerical solution obtained for component 2, see Figure \ref{comp3disperso6} (c), which shows the typical staircase structure of the isotachic train, already seen in the approximate solutions obtained in Experiment 1. The decrease in the quantity of displacer injected prevents the formation of a rectangular pulse for component 1.

\begin{figure}
\includegraphics[width=\textwidth]{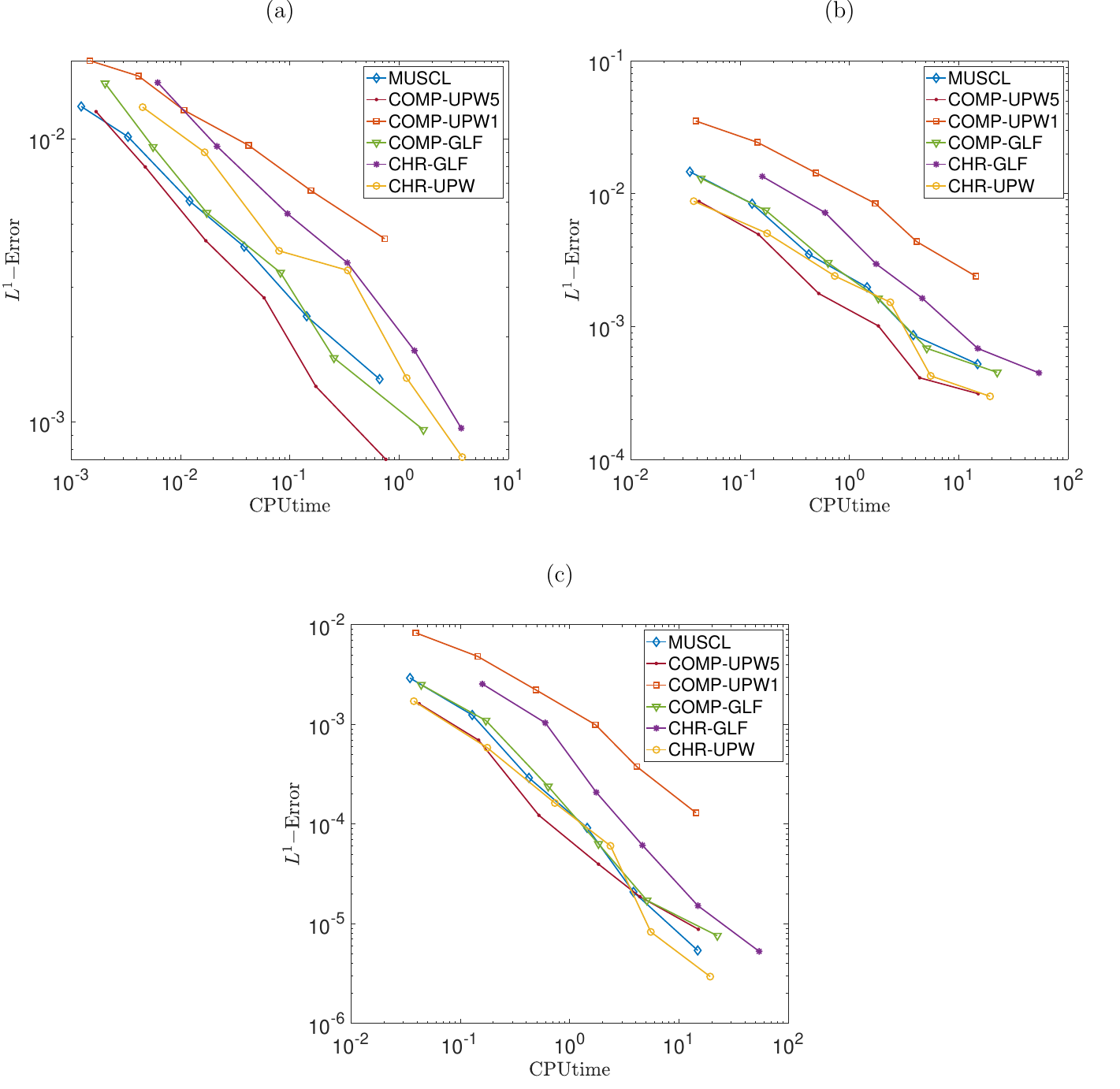}
\caption{Experiment 2. Performance of MUSCL, COMP-UPW5, COMP-UPW1, COMP-GLF, CHR-GLF and CHR-UPW methods for $D_a=0$, $T=1$ (a) and $T=16$ (b). In (c) the performance of the methods with the $L^1$-error computed discarding the $2\%$ of the largest errors for $T=16$ is shown.}
\label{comp3disperso7}
\end{figure}

In Figure \ref{comp3disperso7} we repeat the efficiency analysis performed for the previous experiment, using $\nu=0.9$. The conclusions we extract from it are the same as in Experiment 1. If we discard the largest errors in the computation of the approximate $L^1$-errors, the performance of the CHR-UPW scheme is better than the performance of the other schemes, as the mesh is refined. However, the difference between the errors of all the methods considered in this work is smaller than in Experiment 1, as shown in Figure \ref{comp3disperso5}.  This is probably because the solutions do not develop the isotachic train, therefore the discontinuities are well separated, a fact that contributes to the reduction of the oscillations for the high-order reconstructions.

\subsection{Experiment 3}

\begin{figure}
\includegraphics[width=\linewidth]{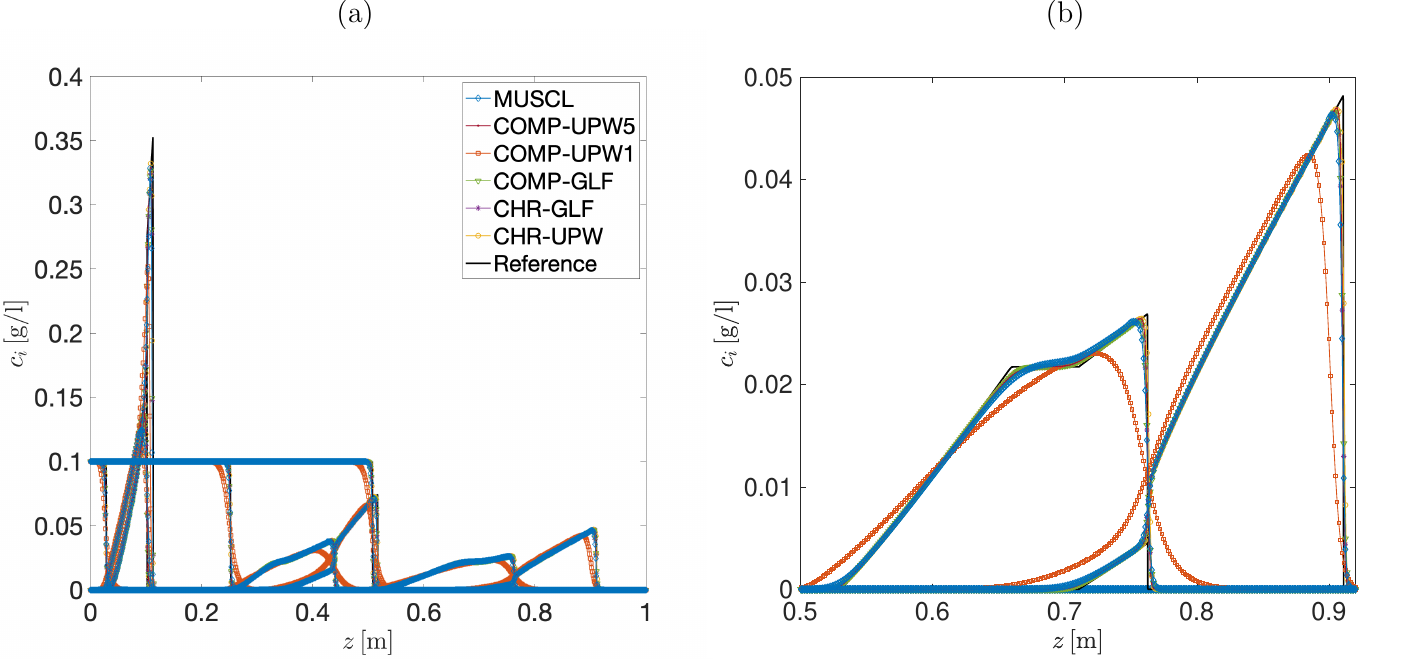}
\caption{Experiment 3. (a) Numerical solutions obtained with MUSCL, COMP-UPW5, COMP-UPW1, COMP-GLF, CHR-GLF and CHR-UPW schemes with $D_a=0$ and $\nu=0.9$ for  $T=1, 8$ and $16$. (b) Enlarged view of the numerical solution of components 1 and 2 in (a) for  $T=16$.}
\label{comp3disperso8}
\end{figure}

\begin{figure}
\includegraphics[width=\textwidth]{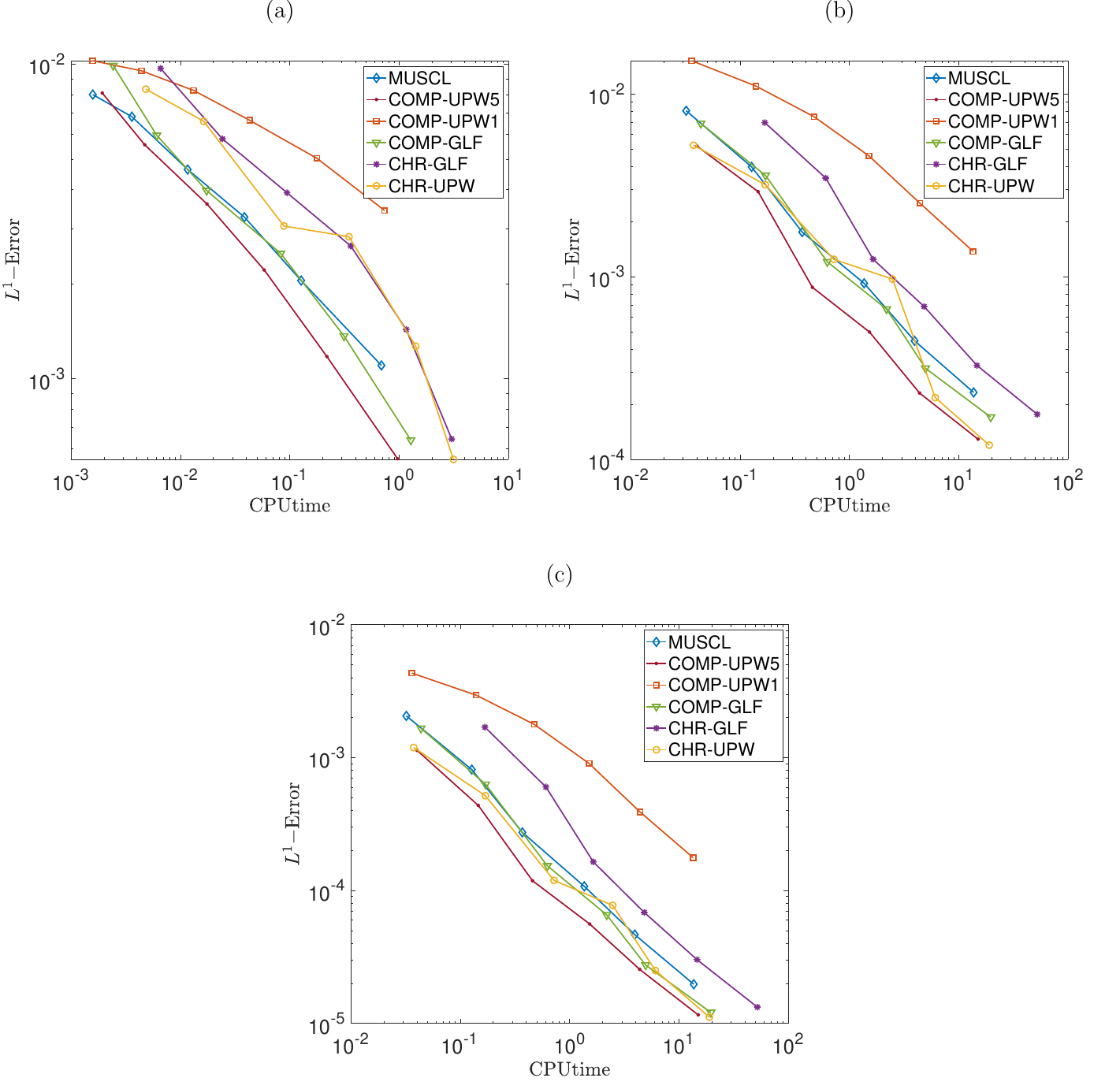}
\caption{Experiment 3. Performance of MUSCL, COMP-UPW5, COMP-UPW1, COMP-GLF, CHR-GLF and CHR-UPW methods for $D_a=0$, $T=1$ (a) and $T=16$ (b). In (c) the performance of the methods with the $L^1$-error computed discarding the $2\%$ of the largest errors for $T=16$ is shown.}
\label{comp3disperso9}
\end{figure}

If we keep reducing the quantity of displacer injected into the column, for instance, we choose $c_3=0.1\,\mathrm{g/l}$, we obtain the numerical solutions in Figure \ref{comp3disperso8}, for $\nu=0.9$. As can be seen, components 1 and 2 do not separate completely in this scenario. Moreover, their profiles do not show rectangular pulses, as is clearly shown by the enlarged view of the numerical solutions for $T=16$ in Figure \ref{comp3disperso8} (b). In this scenario with no isotachic train, the numerical solutions do not show spurious oscillations. The approximate solution obtained with the COMP-UPW1 scheme presents a smoothed-out profile, as was expected and happened also in the previous examples, but the other approximate solutions obtained with the remaining schemes are quite similar, as can be seen in Figure \ref{comp3disperso8} and it is also reflected in the efficiency analysis performed in Figure \ref{comp3disperso9}, with the MUSCL scheme the one with larger errors. In this case, the boost in the performance of the CHR-UPW scheme when discarding the $2\%$ of the largest errors in the computation of the $L^1$-error that could be seen in Figures \ref{comp3disperso5} and \ref{comp3disperso7} does not happen, probably due to the lack of spurious oscillations in the numerical solutions. Moreover, the results obtained with the CHR-UPW scheme in terms of error and computational time are comparable with the ones obtained by the COMP-UPW5 and COMP-GLF schemes.

\subsection{Experiment 4}
\begin{figure}[htb]
\begin{center}
\includegraphics[width=0.6\textwidth]{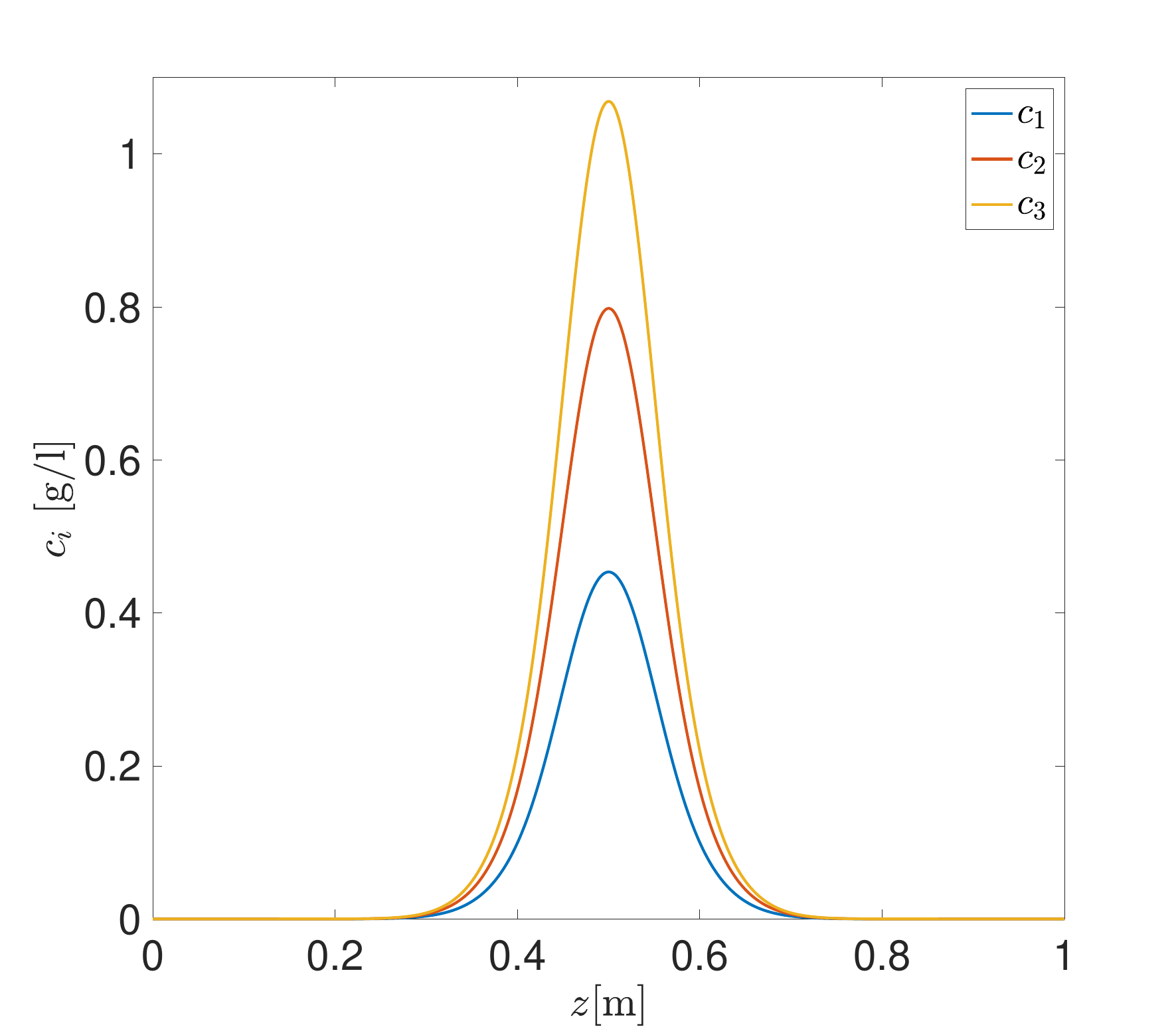}                                                     
\end{center}    
\caption{Experiment 4. Smooth initial data.}
\label{fig:fig_initialcond}
\end{figure}

This experiment,  extracted from \cite{DGM18}, seeks to verify that the CHR-UPW scheme achieves second-order accuracy in simulations featuring smooth solutions. Therefore, we specify smooth initial conditions
\begin{align*}
w_{i}(x, 0)=\rho_{i}\exp(-100(x-1/2)^2),\, i=1,\,2,\,3,
\end{align*}
where $\rho_1=1, \rho_2=2,\rho_3=3$, see Figure \ref{fig:fig_initialcond}. The parameters in the adsorption isotherms are set to $a_1=4, a_2=5, a_3=6$, as in the previous experiments, and $b_i=1,\, i=1,\,2\,,3$. We consider no injection, i.e, $\boldsymbol{c}_{\text{inj}}=0$ in \eqref{eq:bc2w}, and $u=0.2$.

\begin{figure}[htb]
\begin{center}
\includegraphics[width=0.99\textwidth]{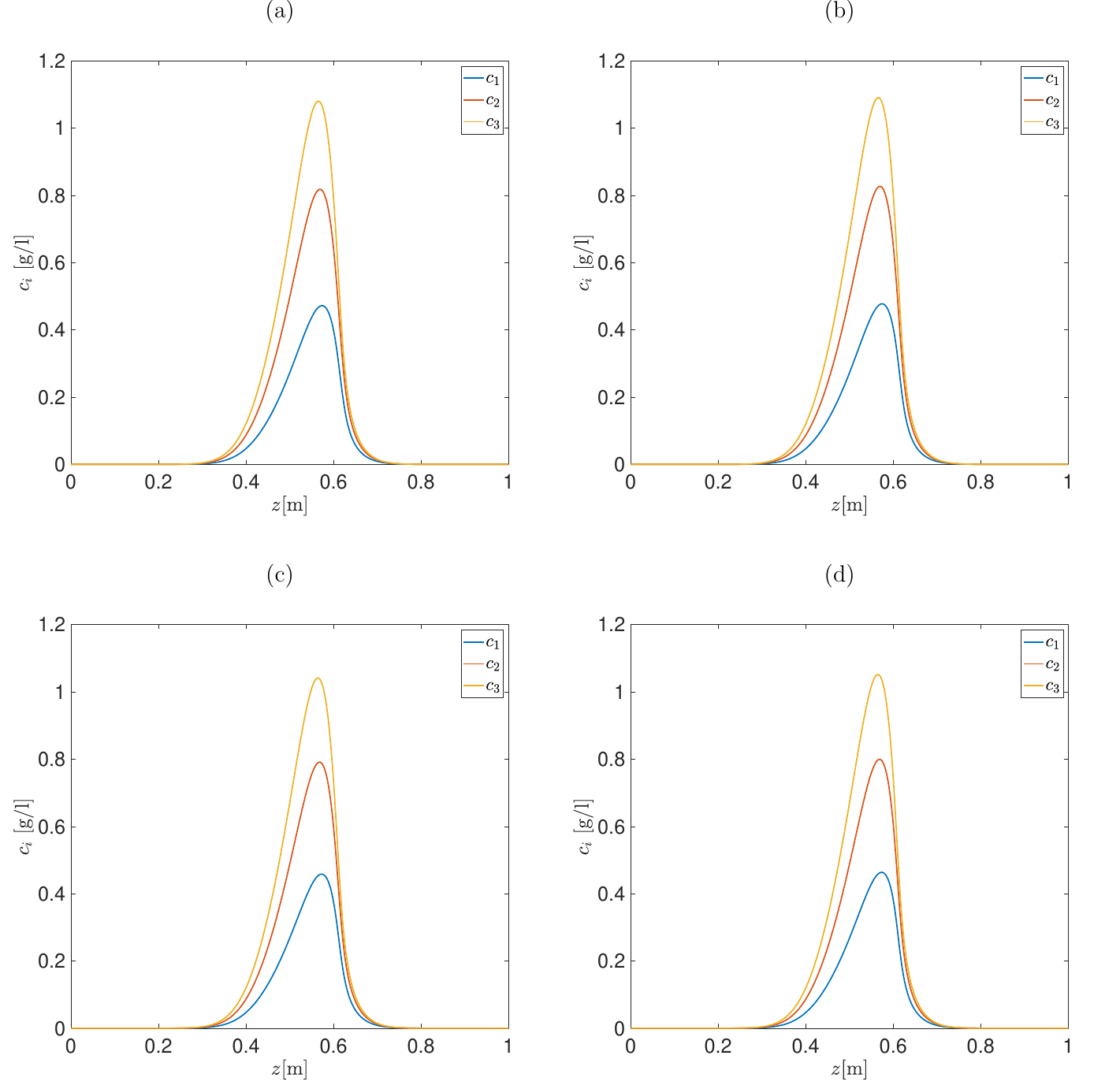}                                                        
\end{center}    
\caption{Experiment 4. Reference solutions of the smooth data test for $T=0.5$ and (a) $D_a=10^{-4}, \nu=0.95$, (b) $D_a=10^{-4}, \nu=1$, (c) $D_a=10^{-5}, \nu=0.95$, (d) $D_a=10^{-5}, \nu=1$.}
\label{fig:fig_order}
\end{figure}

We simulate a short time, until $T=0.5$, before discontinuities appear. In Figure \ref{fig:fig_order}, the reference solutions, obtained with COMP-UPW5 and $m_{\text{ref}}=25600$ cells, are shown for values of $D_a=10^{-4},\,10^{-5}$ and $\nu=0.95,\,1$. We observe in Figure \ref{fig:fig_order} the smooth modified configuration from the initial symmetric configuration, about $z=\mig$, of the variables in Figure \ref{fig:fig_initialcond}.

Based on the approximate errors defined by~\eqref{formula_error} and the mentioned reference solutions, we calculate the numerical order of convergence from pairs of total approximate $L^1$-errors $\smash{e_{m}(T)}$ and $\smash{e_{2m}(T)}$ by 
\begin{align*}    
\theta_m(T) = \log_2 \bigl( e_{m}(T)/ e_{2m}(T) \bigr). 
\end{align*}

In Table \ref{tbl:order_conv}, we show the approximate $L^1$-errors and orders of accuracy for CHR-UPW, COMP-UPW5 and COMP-GLF schemes, for values of $D_a=10^{-4},\,10^{-5}$ and $\nu=0.95,\,1$. We observe that the accuracy orders of all schemes assume values close to two as $m$ increases, as was expected. It is worth mentioning that for this experiment, where shocks and discontinuities are not present in the approximate solutions, component-wise schemes are more accurate than the characteristic-based ones.

\begin{table}[htb]
\begin{center}
\begin{tabular}{|c|c|c|c|c|c|c|c|c|}
\hline
\hline
\multicolumn{9}{|c|}{CHR-UPW}\\ \hline
\hline
&
\multicolumn{4}{|c}{$D_a=10^{-4}$}&
\multicolumn{4}{|c|}{$D_a=10^{-5}$}\\ \hline
&
\multicolumn{2}{|c}{$\nu=0.95$}&
\multicolumn{2}{|c}{$\nu=1$}&
\multicolumn{2}{|c}{$\nu=0.95$}&
\multicolumn{2}{|c|}{$\nu=1$}\\ \hline
$m$&$e_{m}\times 10^{6}$&$	\theta_m$&  $e_{m}\times 10^{6}$&$	\theta_m$&  $e_{m}\times 10^{6}$&$	\theta_m$&  $e_{m}\times 10^{6}$&$	\theta_m$\\ \hline
100 & 1621.02 & 1.77 & 1570.05 & 1.79 & 1688.97 & 1.74 & 1629.61 & 1.75 \\ \hline
200 & 476.04 & 1.92 & 455.49 & 1.92 & 506.58 & 1.91 & 482.93 & 1.90 \\ \hline
400 & 125.65 & 1.99 & 120.47 & 1.98 & 135.02 & 1.99 & 129.44 & 1.98 \\ \hline
800 & 31.62 & 1.98 & 30.63 & 1.97 & 33.88 & 1.99 & 32.80 & 1.99 \\ \hline
1600 & 8.03 & $-$ & 7.81 & $-$ & 8.51 & $-$ & 8.25 & $-$ \\ \hline
\hline
\multicolumn{9}{|c|}{COMP-UPW5}\\ \hline
\hline
100 & 963.81 & 1.73 & 909.76 & 1.73 & 985.31 & 1.70 & 923.35 & 1.69 \\ \hline
200 & 290.32 & 1.97 & 274.51 & 1.96 & 303.42 & 1.96 & 286.59 & 1.95 \\ \hline
400 & 73.97 & 1.99 & 70.42 & 1.99 & 78.07 & 2.00 & 74.38 & 2.00 \\ \hline
800 & 18.58 & 1.96 & 17.71 & 1.96 & 19.54 & 2.00 & 18.64 & 2.00 \\ \hline
1600 & 4.76 & $-$ & 4.56 & $-$ & 4.89 & $-$ & 4.67 & $-$ \\ \hline
\hline
\multicolumn{9}{|c|}{COMP-GLF}\\ \hline
\hline
100 & 932.22 & 1.71 & 842.30 & 1.65 & 941.96 & 1.68 & 855.32 & 1.62 \\ \hline
200 & 284.69 & 1.95 & 268.94 & 1.94 & 294.10 & 1.92 & 278.91 & 1.92 \\ \hline
400 & 73.72 & 1.99 & 70.12 & 1.99 & 77.58 & 1.99 & 73.86 & 1.99 \\ \hline
800 & 18.56 & 1.96 & 17.69 & 1.96 & 19.51 & 2.00 & 18.61 & 2.00 \\ \hline
1600 & 4.76 & $-$ & 4.55 & $-$ & 4.88 & $-$ & 4.66 & $-$ \\ \hline
\end{tabular}
\end{center}
\caption{Experiment 4. Approximate errors $e_m(T)$ for the CHR-UPW, COMP-UPW and COMP-GLF schemes and corresponding numerical orders of convergence $\theta_m$ for the smooth data test with $T=0.5$. The column corresponding to $e_{m}$ has been multiplied by $10^{6}$ to account for the scale.}
\label{tbl:order_conv}
\end{table}

\section{Conclusions} \label{sec:conc}

The present work extends the conservative formulation of the
Equilibrium Dispersive Model developed by Donat, Guerrero and Mulet in
\cite{DGM18} for {\it Langmuir} adsorption isotherms, to the generalized
{\it Langmuir}-type adsorption isotherms that we propose. We have proven that for this family of functions, there is a smooth bijection between the
concentrations of the solutes in the liquid phase and the conserved variables that
allows us to write the model as a well-posed system of conservation laws with
diffusive corrections. This correspondence allows us to numerically recover the characteristic information of the Jacobian matrix of the convective fluxes and use it to design an implicit-explicit scheme that uses characteristic-based numerical fluxes for the fifth-order WENO reconstruction technique.

As reflected by the numerical examples performed using T\'oth's isotherms, the use of  {cha\-rac\-te\-ris\-tic-based} numerical fluxes is an excellent option to eliminate spurious oscillations caused by the reconstruction method, although the computational time needed to obtain the approximate solutions may be higher. The numerical results displayed allow us to conclude that the proposed numerical technique is a reliable and robust tool for numerically solving this model.

In this paper, we have applied the proposed scheme to the equilibrium
dispersive model for a fixed bed. However,  it might be adapted to
more general models of chromatographic columns and also to simulating
moving beds (SMB) models.   Future work will also explore
adsorption isotherms beyond T\'oth's isotherms.

 \section*{Acknowledgments}

This research has been partially supported by grant PID2020-117211GB-I00 funded by Ministerio de Ciencia e Investigaci\'on MCIN/AEI/10.13039/501100011033 and by Conselleria de Innovaci\'on, Universidades, Ciencia y Sociedad Digital through project CIAICO/2021/227.

We are grateful to Francisco Guerrero Cortina for preliminary work on
this project.

\end{document}